
\input amstex

\magnification 1200
\loadmsbm
\parindent 0 cm

\define\nl{\bigskip\item{}}
\define\snl{\smallskip\item{}}
\define\inspr #1{\parindent=20pt\bigskip\bf\item{#1}}
\define\iinspr #1{\parindent=27pt\bigskip\bf\item{#1}}
\define\einspr{\parindent=0cm\bigskip}

\define\co{\Delta}
\define\st{$^*$-}
\define\ot{\otimes}

\define\C{\Bbb C}

\input amssym

\centerline{\bf Compact and discrete subgroups of algebraic quantum groups I}
\bigskip

\centerline{\it M.\ B.\ Landstad  \rm($^*$) and \it A.\ Van Daele
\rm ($^{**}$)}
\bigskip

{\bf Abstract} \nl Let $G$ be a locally compact group. Consider
the C\st algebra $C_0(G)$ of continuous complex functions on $G$,
tending to $0$ at infinity. The product in $G$ gives rise to a
coproduct $\co_G$ on the C\st algebra $C_0(G)$. A locally compact
{\it quantum} group is a pair $(A,\co)$ of a C\st algebra $A$ with
a coproduct $\co$ on $A$, satisfying certain conditions. The
definition guarantees that the pair $(C_0(G),\co_G)$ is a locally
compact quantum group and that conversely, every locally compact
quantum group $(A,\co)$ is of this form when the underlying C\st
algebra $A$ is abelian.
\smallskip

Assume now that $G$ is a locally compact group with a compact open
subgroup $K$. A continuous complex function $f$ with compact
support on $G$ is called of {\it polynomial type} if there
exist finitely many continuous functions $(f_i)$ on $G$ and
$(g_i)$ on $K$ such that $f(pk)=\sum f_i(p)g_i(k)$ for all $p\in
G$ and $k\in K$. The set $P(G)$ of such functions is a dense \st
subalgebra of $C_0(G)$, independent of the choice of $K$. The
comultiplication $\co_G$ leaves $P(G)$ invariant and the pair
$(P(G),\co_G)$ is a so-called {\it algebraic quantum group}.
Algebraic quantum groups are a special class of locally compact
quantum groups. This class contains the discrete quantum groups,
the compact quantum groups and it is self-dual. As mentioned
above, also a locally compact group with a compact open subgroup
falls into this class. In fact, if a locally compact group $G$
gives rise to an algebraic quantum group, it must contain a
compact open subgroup.
\smallskip

When $K$ is a compact open subgroup of $G$, the characteristic
function $\chi_K$ of $K$ is a self-adjoint idempotent $h$
satisfying $\co_G(h)(1\ot h)=h\ot h$. In this paper, we consider
such {\it group-like projections} in a general algebraic quantum
group. In other words, we study compact quantum subgroups of
algebraic quantum groups. For a general algebraic quantum group,
we look at the analogues of the \st algebra of continuous complex
functions on the subgroup $K$ and the \st algebra of complex
functions with finite support on the quotient space $G/K$ in the
group case, as well as their dual objects. In a forthcoming paper on 
this subject, we plan to include more examples and special cases to illustrate
the different notions and results of this paper.
\smallskip

Results in algebraic quantum groups usually are also true (in an
adapted form) in general locally compact quantum groups and we
discuss further research in this direction. In particular, we will
propose a definition of a {\it totally disconnected (locally
compact) quantum group}. \bigskip

April 2007 
\bigskip

\hrule
\medskip
($^*$) Department of Mathemactics and Statistics, University of Trondheim, N-7491
Trondheim (Norway). E-mail: Magnus.Landstad\@math.ntnu.no
\smallskip
($^{**}$) Department of Mathematics, K.U.\ Leuven, Celestijnenlaan 200B,
B-3001 Heverlee (Belgium). E-mail: Alfons.VanDaele\@wis.kuleuven.be

\newpage

\bf 0. Introduction \rm
\nl
This paper studies some material related to the theory of
locally compact quantum groups. Let us therefore first recall this
concept.
\nl
\it Locally compact quantum groups \rm
\nl
Let $G$ be a locally compact group. Consider the C\st algebra
$C_0(G)$ of continuous complex functions on $G$ tending to $0$ at
infinity. It is well-known that the topological structure of $G$
is encoded in the C\st algebraic structure of $C_0(G)$. The
product in $G$ gives rise to a coproduct on $C_0(G)$ in the
following way. First we identify the minimal C\st tensor product
$C_0(G) \ot C_0(G)$ of $C_0(G)$ with itself with the C\st algebra
$C_0(G \times G)$ of continuous complex functions, tending to $0$
at infinity on the cartesian product $G \times G$ of $G$ with
itself. Then we identify the C\st algebra $C_b(G\times G)$ of
continuous bounded complex functions on $G\times G$ with the
multiplier algebra $M(C_0(G) \otimes C_0(G))$ of $C_0(G) \otimes
C_0(G)$. Now we define a \st homomorphism $\co : C_0(G) \to
M(C_0(G) \ot C_0(G))$ by $(\co(f))(p,q) = f(pq)$ whenever $f\in
C_0(G)$ and $p,q\in G$. This \st homomorphism is a
comultiplication on the C\st algebra $C_0(G)$ and it carries all
the information about the product in $G$. Indeed, the product in
$G$ can be recovered from this coproduct on $C_0(G)$.
\snl
The theory of locally compact quantum groups is motivated by the
above observation. Now a pair $(A,\co)$ of a C\st algebra $A$ (not
necessarily abelian) and a comultiplication $\co$ on $A$ is
considered. The comultiplication $\co$ on $A$ is a \st
homomorphism from $A$ to the multiplier algebra $M(A \ot A)$ of
the minimal C\st tensor product $A\ot A$ satifsying certain
properties (such as coassociativity).
\snl
The pair $(A,\co)$ will be a locally compact quantum group
provided there exist left and right Haar weights. These are the
non-commutative analogues of the left and the right Haar measures
in the group case. In the general quantum case, the existence of
such invariant weights is a part of the definition.
\snl
For more information about the general theory of locally compact
quantum groups, we refer to [K-V1], [K-V2] and [K-V3]. See also [VD6] and [VD7] 
for a new, von Neumann algebraic approach to this theory.

\nl
\it About this paper \rm

\nl
Very little knowledge about general locally compact quantum groups however is assumed
for understanding this paper. Some notions of this
general theory will help to have a better idea of the main results
of this work. This paper should be readable (and of interest) for algebraists, familiar with some Hopf algebra theory, as well as operator algebraists, with some knowlegde about the operator algebra approach to quantum groups.

\snl
Indeed, this paper mainly deals with multiplier Hopf \st algebras
with positive integrals; the so-called {\it algebraic quantum
groups} (cf.\ [VD3]). It is well-known that the purely algebraic aspects of
locally compact quantum groups can be studied - up to a certain
level - within the special case of such algebraic quantum groups
(see [VD5]). In fact, historically, the development of the
general theory was greatly motivated by the results on algebraic
quantum groups. In the same spirit, the material studied in this
paper is also expected to contribute to the further development of
the general theory of locally compact quantum groups.

\nl
\it Motivation \rm

\nl
The interest for what we are doing here, finds its origin in a
question related with the above discussed role of the algebraic
quantum groups in the general theory. Let $G$ be a locally compact
group. When is the associated pair $(C_0(G),\co)$ an algebraic
quantum group? More precisely, when is there a dense \st
subalgebra $A$ of $C_0(G)$ such that $(A,\co)$ is an algebraic
quantum group? This is the case when $G$ is discrete. Then $A$ is
the \st algebra of complex functions on $G$ with finite support.
It is also the case when $G$ is a compact group. Now $A$ is the
\st algebra of polynomial functions (i.e.\ matrix elements for
finite-dimensional representations, considered as complex
functions on the group $G$). And of course, trivial cases build
from these two special cases, like the direct product of a
discrete and a compact group, will also give rise to an algebraic
quantum group.

\snl
In [L-VD1] an answer to this question is given. It turns out
that $C_0(G)$ is an algebraic quantum group (in the above sense),
if and only if the group $G$ contains a compact open subgroup.

\snl
Now, if $K$ is a compact open subgroup of a locally compact group
$G$, and if we let $h$ be the characteristic function $\chi_K$ of
$K$, then $h$ is a projection (i.e.\ a self-adjoint idempotent) in
$C_0(G)$ and it satisfies the equation $\co(h)(1\ot h)=h\ot h$.
One of the main objects that we study in this paper  are precisely
such projections in multiplier Hopf \st algebras.

\snl
Throughout the paper, we will have this motivating example in mind. We 
will use it to illustrate the various objects that we introduce and 
the properties that we prove. In two forthcoming papers, we treat some 
other, more complicated examples together with some special cases 
(see [L-VD2] and [L-VD3]).

\nl
Later in this introduction, we will recall some notions and
conventions, used in this paper, together with the basic
references. We also refer to the appendix for a review of the
basic properties of multiplier Hopf algebras.

\nl
\it Content of the paper\rm

\nl
In {\it Section 1} of this paper, we will initiate the study of
what we will call {\it group-like} projections in a multiplier Hopf
\st algebra. These are (non-zero) elements $h$ in a multiplier
Hopf \st algebra $A$ satisfying $h^2=h=h^*$ (self-adjoint
idempotents or projections) and moreover $\co(h)(1\ot h)=h\ot h$.
We call such a projection group-like because, in the abelian case,
it really comes from a subgroup (see Proposition 1.4 in Section
1). We found it quite remarkable how many properties of such a
group-like projection can be proven from these elementary
assumptions. This is done in Section 1.

\snl
In {\it Section 2}, we begin with showing that, if $h$ is a
{\it central} group-like projection in $A$ (that is if $ha=ah$ for all
$a\in A$), then, as expected, the \st algebra $Ah$ can be made
into a compact quantum group by cutting down the multiplication,
i.e.\ by defining $\co_0:Ah \to Ah\ot Ah$ as $\co_0(a)=\co(a)(h\ot
h)$ when $a\in Ah$. The integral on $Ah$ is obtained by
restricting the integrals of $A$. This result is not very
difficult. The rest of the section is devoted to the general case,
where $h$ is no longer assumed to be central. Now we take $hAh$ as the
underlying algebra and we have to cut down $\co$ on both sides
with $h\ot h$. We are left with a positive
comultiplication which is no longer a homomorphism. We get a
so-called {\it compact quantum hypergroup} (or an algebraic quantum hypergroup of compact type in the sense of [De-VD1]). 

\snl
In {\it Section 3} we study the equivalent of the algebra of
functions that are constant on right cosets. More precisely, we
study the 'right leg' of $\co(h)$ for a group-like projection $h$. 
It can be characterized as the
smallest subspace $C$ of $A$ such that $\co(h)(A\ot 1)\subseteq
A\ot C$. It turns out that $C$ is a \st subalgebra of $A$. It is
also left invariant in the sense that the right leg of $\co(C)$ is
contained in $C$. In the event that $C$ is invariant under the
antipode (this is the case when the left and the right leg of
$\co(h)$ are the same), then $(C,\co)$ is a discrete quantum
group. In the abelian case, this happens when the compact open
subgroup $K$ is normal. Then $C$ is the discrete quantum group of
functions of finite support on the quotient group $G/K$, now a
discrete group. Again, in the general case, we  have to cut down
the comultiplication, still restrict to a smaller \st subalgebra
$C_1$ and then we get what can be considered as a {\it discrete quantum
hypergroup} (or an algebraic quantum hypergroup of discrete type in the 
sense of [De-VD1]).

\snl
Throughout these three preceding sections, we also consider the case of a pair of
group-like projections, one smaller than the other. This is important for later work on
what will be defined as a totally disconnected locally compact quantum group.
\snl
In the last section, {\it Section 4}, we formulate some
conclusions and we discuss future research.
\snl
In an {\it Appendix} we first collect the basic definitions and
properties about multiplier Hopf $^*$-algebras $(A,\Delta)$, algebraic quantum
groups (multiplier Hopf $^*$-algebras with positive integrals) and
duality. We also introduce and study the notion of the legs of
$\Delta(a)$ for elements $a\in A$. We define left invariant subalgebras and prove
some properties. These objects and results are used in the paper,
but as they are of some general interest and not dependent upon
the notion of a group-like projection, we have put this material in
the appendix.

\nl
\it Notions and conventions, basic references \rm

\nl
As we have mentioned already, this paper will mainly deal with
purely algebraic objects. We will work with algebras $A$ over
$\Bbb C$, possibly without identity, but always with a
non-degenerate product.  We use $A'$ for
the space of linear functionals on $A$. We will use $M(A)$ to denote the
multiplier algebra of $A$. The identity in an algebra, e.g.\ in $M(A)$, will be denoted by $1$ (while
$e$ will be used to denote the identity in a group). The identity map, say 
from $A$ to itself, will always be denoted by $\iota$. The (algebraic) tensor product $A\otimes A$ is 
again an algebra over $\Bbb C$ with a non-degenerate product. We identify e.g.
$A\otimes \Bbb C$ with $A$ and therefore, the slice map $\iota\otimes\omega$ for an element
$\omega\in A'$ is considered to be a map from $A\otimes A$ to $A$. 

\snl
A {\it comultiplication} on $A$ is a
homomorphism $\co:A \to M(A \ot A)$ satisfying certain properties (such as coassociativity).
A {\it multiplier Hopf (\st)algebra} is a (\st)algebra $A$ with a
non-degenerate product and a comultiplication $\co$ satisfying
certain conditions. A multiplier Hopf algebra is called {\it
regular} if the opposite comultiplication $\co^{\text{op}}$,
obtained by composing $\co$ with the flip, also satisfies these
properties. In the case of a \st algebra it is assumed that $\co$ is a $^*$-homomorphism and regularity is automatic.
For a multiplier Hopf algebra $(A,\co)$, we have the existence of
a unique {\it counit} and a unique {\it antipode}.

\snl
Any Hopf (\st)algebra is a multiplier Hopf (\st)algebra.
Conversely, if $(A,\co)$ is a multiplier Hopf (\st)algebra and if
$A$ has an identity, it is a Hopf (\st)algebra. So we see that the
theory of multiplier Hopf algebras extends in a natural way the
theory of Hopf algebras to the case where the underlying algebras
are not required to have an identity. The use of Sweedler's notation in the case of Hopf algebras is a common practice. It is also justified in the case of multiplier Hopf algebras (see e.g.\ [Dr-VD] and [Dr-VD-Z]) and whenever convenient, we will also do so in this paper.

\snl
Let $(A,\co)$ be a regular multiplier Hopf algebra. A linear
functional $\varphi$ on $A$ is called left invariant if
$(\iota\ot\varphi)\co(a)=\varphi(a)1$ for all $a\in A$.  A {\it
left integral} is a non-zero left invariant linear functional on
$A$. Similarly, a non-zero right invariant linear functional on
$A$ is called a {\it right integral}.

\snl
If $(A,\co)$ is a multiplier Hopf \st algebra with a positive left
integral, then there is also a positive right integral. We will
use the term {\it algebraic quantum group} for such a multiplier
Hopf \st algebra. In this paper, we mainly deal with algebraic
quantum groups, defined in this sense. However, although for many of the 
results we need an integral, it is not so clear if the positivity of
this integral is important. Still, we make the convention that integrals
are assumed to be positive. This implies that our $^*$-algebras are in
fact operator algebras, i.e.\ that they can be represented as $^*$-algebras of bounded operators on a Hilbert space.

\snl
We have collected the most important definitions and properties
about multiplier Hopf $^*$-algebras and algebraic quantum groups
in an appendix. For more details, we refer to [VD1], [VD3] and
[VD-Z]. For results about ordinary Hopf algebras, we refer to the basic works [A] and [S].

\snl
We will freely use results from these basic references in this
paper. In particular, we will use $\varphi$ to denote a left
integral and we use $\psi$ for a right integral. We use $\sigma$
for the {\it modular automorphism} of $\varphi$, satisfying
$\varphi(ab)=\varphi(b\sigma(a))$ for all $a,b\in A$. Similarly
$\sigma'$ is used for the modular automorphism of $\psi$. The {\it
modular element} $\delta$ is a multiplier in $M(A)$, defined and
characterized by $(\varphi\ot\iota)\co(a)=\varphi(a)\delta$ for
all $a\in A$. It is invertible and the inverse satisfies
$(\iota\ot\psi)\co(a)=\psi(a)\delta^{-1}$. Finally, there is the
{\it scaling constant} $\nu$ defined by
$\varphi(S^2(a))=\nu\varphi(a)$ where $S$ is the antipode. We also have
collected some of the main formulas relating these various
objects, associated with the original pair $(A,\Delta)$ as well as 
with the dual pair $(\hat A,\hat \Delta)$, in the appendix.

\nl\nl

\bf Acknowledgements \rm
\medskip

This research has been going on over a long period, starting at The Centre for Advanced Study (CAS) in Oslo, and continuing at the K.U.Leuven and NTNU. We thank for their hospitality and financial support. In addition we also received support from The Research Council of Norway (NFR), the FWO-Vlaanderen and The Research Council of the K.U.Leuven.
The second named author (Alfons Van Daele) is also grateful to his
colleague Johan Quaegebeur at the University of Leuven for taking
care of most of his teaching when he was on sabbatical (2002-2003).
\nl\nl

\bf 1. Group-like projections in multiplier Hopf \st algebras \rm
\nl
Let $(A,\co)$ be a multiplier Hopf \st algebra. We begin with the
definition of a group-like projection. This is essentially the
basic object for the rest of the paper. We illustrate the notion by 
some basic examples and a simple, but important result (Proposition 1.4 below).
We prove elementary properties of such group-like projections. We also show
that, in the case of algebraic quantum groups, the 'Fourier transform' of a 
group-like projection is again a group-like projection in the dual $(\hat A,
\hat\Delta)$.
\nl
\it Definition and examples of group-like projections \rm

\inspr{1.1} Definition \rm
Let $h$ be a self-adjoint projection in $A$. We call $h$ {\it
group-like} if
$$ \co(h)(1 \ot h)= h \ot h. $$

\einspr
That $h$ is a self-adjoint projection means that $h^2=h=h^*$. 
Recall also that in a multiplier Hopf algebra, we have $\co(a)(1 \ot
b)\in A \ot A$ for all $a,b\in A\ot A$ so that the equation in the
definition in principle can occur.
\snl
We will always assume that $h\neq 0$ when we consider a group-like
projection $h$ in this paper.
\snl
Recall that an element $a$ in a Hopf algebra is called group-like
if $\co(a)=a\ot a$. In general, this notion does not make much
sense in a multiplier Hopf algebra. One rather has to look at
elements in the multiplier algebra $M(A)$ of $A$ with this
property. Here however, we use the term in a different sense. We
call the element $h$ group-like because, in the classical case, it
behaves like the characteristic function of a subgroup. This was
mentioned already in the introduction. It will also become clear
from the following examples. They will motivate and justify the
use of this terminology.
\snl
The first two examples are very simple and easy to understand. We
include them mainly for didactical purposes.

\inspr{1.2} Examples \rm i) Let $G$ be a (discrete) group and let
$A$ be the algebra of complex functions on $G$ with finite support
with pointwise operations. Consider the obvious comultiplication
$\co$ on $A$ given by $(\co(f))(p,q)=f(pq)$ when $f\in A$ and $p,q
\in G$. Let $K$ be any finite subgroup and let $h$ be the
characteristic function on $K$, i.e.\ $h(p)=1$ if $p\in K$ and
$h(p)=0$ if $p \notin K$. Then $(\co (h)(1 \ot
h))(p,q)=h(pq)h(q)$. This will be non-zero if both $pq$ and $q$
are in $K$. This is true if and only if both $p$ and $q$ are in
$K$. So $\co (h)(1 \ot h)=h \ot h$.
\snl ii) Again let $G$ be a discrete group, but now consider the
group algebra $B$ over $\C$. The comultiplication $\co$ is given
by $\co(\lambda_p) = \lambda_p \ot \lambda_p$ where $p \mapsto
\lambda_p$ is the canonical imbedding of $G$ in $B$. Again let $K$
be a finite subgroup of $G$. Now put
$$ h = \frac{1}{n}\sum_{p\in K}\lambda_p$$
where $n$ is the number of elements in $K$. It is easy to verify
that $h$ is a self-adjoint projection (observe that
$\lambda_p^*=\lambda_{p^{-1}}$) and that it satisfies $\co(h)(1
\ot h) =h\ot h$.
\einspr

Here are the (more interesting) topological equivalents of the
examples in 1.2.

\inspr{1.3} Examples \rm i) Let $G$ be a locally compact group.
Consider the C\st algebra $C_0(G)$ of continuous complex functions
on $G$ vanishing at infiniy. The comultiplication $\co$ is now
defined by $(\co(f))(p,q)=f(p,q)$ when $f \in C_0(G)$ and $p,q \in
G$. This will not be a multiplier Hopf \st algebra in general.
However, if there is a compact open subgroup $K$ of $G$, then it
is shown in [L-VD1] that there is a dense \st subalgebra $A$ of
$C_0(G)$ such that $(A,\co)$ is a multiplier Hopf \st algebra with integrals 
(i.e.\ an algebraic quantum group).
\snl 
The characteristic function $h$ of $K$ will be a self-adjoint
projection in $A$ satisfying the equation $\co(h)(1 \ot h)=h \ot
h$. 
\snl 
ii) Let $G$ be as in i), that is a locally compact group
with a  compact open subgroup $K$. Now consider the reduced group
C\st algebra $C_r^*(G)$. This is the C\st algebra generated by
operators of the form
$$\int_{p \in G}f(p)\lambda_p dp$$
where integration is over the left Haar measure on $G$, where
$\lambda_p$ is the translation operator on $L^2(G)$, defined by
$(\lambda_p\xi)(q)=\xi(p^{-1}q)$ when $\xi \in L^2(G)$ and $p,q \in G$
and where $f \in  C_c(G)$, the space of continuous complex functions
with compact support on $G$. The comultiplication on $C_r^*(G)$ is given
by
$$\co\left(\int f(p) \lambda_p dp\right)=
\int f(p)(\lambda_p \ot \lambda_p) dp$$
for $f \in C_c(G)$. Also in this case, there exists a dense $^*$-subalgebra that
is a multiplier Hopf $^*$-algebra with integrals (see [L-VD1]).
\snl
Again $h$, defined by
$$h=\frac{1}{n}\int_{p\in K}\lambda_p dp,$$
where $n$ is the Haar measure of $K$, will be a group-like projection
in the sense of Definition 1.1.
\einspr

In both cases, the dense $^*$-subalgebra is constructed by taking \it polynomial 
functions \rm on the group $G$ w.r.t.\ the open compact subgroup $K$. These are continuous 
complex functions $f$ on $G$ with compact support such that there exist finitely many 
continuous functions $(f_i)$ on $G$ and $(g_i)$ on $K$ satisfying
$$f(pk)=\sum_i f_i(p)g_i(k)$$
for all $p\in G$ and $k\in K$. Again see [L-VD1].
\nl
There is also the following result, related with the first example above.
Because it is rather fundamental for the motivation of the basic definition (Definition 1.1 of a group-like projection) and hence for the rest of this  paper, let us formulate it as a separate result.

\inspr{1.4} Proposition \rm
Let $G$ be a locally compact group and
let $K$ be a compact and open subset of $G$. Let $h$ be the
characteristic function of $K$. Then $K$ is a subgroup if and only
if $h$ is a group-like projection.

\snl\bf Proof: \rm
First assume that $\co(h)(1 \ot h) = h \ot h$. Take points $p,q$ in
$G$. If $pq \in K$ and $q \in K$, then
$$h(p)=h(p)h(q)=(h\ot h)(p,q)=(\co(h)(1\ot h))(p,q)=h(pq)h(q)=1.$$
Hence $p\in K$. So $KK^{-1}\subseteq K$. As it is assumed that $h\neq
0$, we have $K$ non-empty and so it follows that the identity $e$ of
the group must be in $K$.  Then it follows that $K^{-1}\subseteq
K$. Taking inverses, we also get the other inclusion and so $K^{-1}=K$. Combined
with the first result, namely that $KK^{-1}\subseteq K$, we finally
get $KK\subseteq K$ and we see that $K$ is a group.
\snl
Conversely, if $K$ is a group, then for any pair of elements $p,q\in
G$ we have $p,q\in K$ if and only if $pq, q\in K$ and then $\co(h)(1 \ot h)=h \ot
h$.
\einspr
Apart from these motivating examples, we will also keep in mind some
special cases:

\inspr{1.5} Examples \rm i) Assume that $(A,\co)$ is of compact
type (see Definition 5.2 in [VD3]), that is when $A$ has an
identity and hence when $(A,\co)$ is actually a Hopf \st algebra.
Then obviously 1 is a group-like projection in the sense of the
above definition. \snl ii) Assume that $(A,\co)$ is of discrete
type (see again Definition 5.2 in [VD3]), that is when there
exists a non-zero element $h\in A$ satisfying $ah=\varepsilon(a)h$
for all $a\in A$. Then
$$ \co(h)(1\ot h)=(\iota\ot\varepsilon)\co(h) \ot h = h \ot h$$
so that $h$ is group-like.
\einspr

These are two extreme cases. Taking the tensor product of a
compact type with a discrete type, we get an example 'in between'
by considering $1\ot h$. In the group case (Example 1.3.i) this
would correspond to taking a direct product of a compact group
with a discrete group and looking at the compact group as a
subgroup of this direct product.
\snl
Most of these examples are rather trivial. Even the two examples
coming from a locally compact group with a  compact open subgroup
are not very difficult. We will use these examples throughout the
paper for illustrating various results. However, all these
examples are in a way too simple for a complete understanding.
In [L-VD2] we  give more and less trivial examples. One of them is the quantum double. In a 
sense, this is a twisted version of the example mentioned above of
the tensor product of a compact type with a discrete type.

\nl
\it Properties of group-like projections in algebraic quantum
groups \rm

\nl
Now, we give the first properties of such group-like projections.
We will mainly deal with algebraic quantum groups.

\inspr{1.6} Proposition \rm
Let $(A,\co)$ be an algebraic quantum group and $h$ a group-like
projection in $A$. Then also $\co(h)(h \ot 1) = h \ot h$. Moreover
$\varepsilon(h)=1$ and $S(h)=h$ (where $\varepsilon$ is the counit and
$S$ the antipode of $A$).

\snl\bf Proof: \rm
Apply $\varepsilon\ot\iota$ to the defining
equation $\co(h)(1 \ot h) = h \ot h$. Then we get
$h^2=\varepsilon(h)h$. As $h^2=h$ and $h$ is assumed to be
non-zero, we obtain $\varepsilon(h)=1$. \snl Now, apply
$\iota\ot\co$ to the defining equation and use the Sweedler
notation (see the appendix) to obtain
$$(h_{(1)} \ot h_{(2)} \ot h_{(3)})(1\ot \co(h))=h \ot \co(h).$$
Then apply $S$ on the first leg, and multiply the first leg with the
second one. We get
$$\align
       (S(h) \ot 1)\co(h)
             &= (S(h_{(1)})h_{(2)}\ot h_{(3)}) \co(h) \\
             &= (1 \ot h)\co(h)\\
             &= h \ot h.
\endalign$$
The last equality follows from the defining equation by taking
adjoints.
\snl
Now apply $\iota\ot\varphi$, where $\varphi$ is a positive left
integral, to get $S(h)\varphi(h)=h\varphi(h)$. Because $\varphi$ is
faitful and $h\neq 0$, we have $\varphi(h)=\varphi(h^*h)\neq 0$. Therefore $S(h)=h$.
\snl
Combining with the above equation, we now get
$$(h\ot 1)\co(h) = h \ot h$$
and taking adjoints, we get the desired equation.
\einspr

It is not clear if these results will still be true for any
multiplier Hopf \st algebra (not an algebraic quantum group).
Clearly still $\varepsilon(h)=1$ will hold. When we apply
$\iota\ot\varepsilon$ to the equation $(S(h)\ot 1)\co(h)=h\ot h$,
we get $S(h)h=h$. If we take adjoints, we get
$h=hS(h)^*=hS^{-1}(h)$. If we apply again $S$, we get
$S(h)=hS(h)$. If we also assume that $\co(h)(h\ot 1) = h \ot h$,
we get using a similar argument, that $h=hS(h)$. Then also
$S(h)=h$.
\snl
In some sense, the result proven earlier in Proposition 1.4, can be seen as an
illustration of Proposition 1.6. Indeed, using the notations of 1.4 and 1.6, we see that
$\varepsilon(h)=1$ means $e\in K$ while $S(h)=h$ means $K^{-1}=K$.
\nl
The following results relate the various data of an algebraic quantum
group with such a group-like projection.

\inspr{1.7} Proposition \rm
Let $(A,\co)$ be an algebraic quantum group and let $h$ be a group-like
projection in $A$. Then the scaling constant $\nu$ is $1$. We have
$h\delta=h$ where $\delta$ is the modular element relating the left
and the right integrals. Also $\sigma(h)=\sigma'(h)=h$ where $\sigma$
and $\sigma'$ are the modular automorphisms for respectively the left
and the right integrals.

\snl\bf Proof: \rm
Because $S(h)=h$ we get $S^2(h)=h$ and so
$\varphi(h)=\varphi(S^2(h))=\nu\varphi(h)$ where $\varphi$ is the left
integral. We have seen already that $\varphi(h)\neq 0$. So $\nu=1$.
\snl
If we apply $\varphi\ot\iota$ to $\co(h)(1\ot h)=h \ot h$, we get
$\delta h=h$ as $(\varphi\ot\iota)\co(h)=\varphi(h)\delta$. Taking
adjoints gives $h\delta=h$.
\snl
We have $\co(h)=\co(S^2(h))=(\sigma\ot{\sigma'}^{-1})\co(h)$ (see Lemma 3.10 in
[K-VD] and Proposition A.6 of the appendix). If we multiply on the left with $h\ot 1$ and on the right with
$1\ot {\sigma'}^{-1}(h)$ we get
$$h \ot h {\sigma'}^{-1}(h)=h\sigma(h) \ot {\sigma'}^{-1}(h).$$
So $h\sigma(h)$ is a scalar multiple of $h$. However, as
$\varphi(h\sigma(h))=\varphi(hh)=\varphi(h)$, we must have
$h\sigma(h)=h$.
Taking
adjoints, and using that $\sigma(h)^*=\sigma^{-1}(h^*)$ (cf.\ Appendix), we also get $h=\sigma^{-1}(h)h$ and if we apply $\sigma$, we find $\sigma(h)=h\sigma(h)$. Hence, as also $h=h\sigma(h)$, we get
$\sigma(h)=h$. Similarly $\sigma'(h)=h$.
\einspr

We see that the existence of a group-like projection implies that
the scaling constant is trivial. For general locally compact
quantum groups, we know examples with a non-trivial scaling
constant (cf.\ [VD4]). For algebraic quantum groups however, the
question was  open for a long time. There were no known examples of algebraic
quantum groups with non-trivial scaling constants. The
above result is an indication that such examples can not be
found. Indeed, very recently, it is shown in [DC-VD] that for algebraic quantum groups (with positive integrals), the scaling constant is always trivial. Observe however that there are multiplier Hopf algebras with integrals and non-trivial scaling constant (see some of the examples in [VD3]).
\nl
It is not possible to illustrate the previous
proposition with the given examples. These are too simple for this
purpose. Indeed, $\nu=1$ in these examples anyway and the modular
automorphisms $\sigma$ and $\sigma'$ are also trivial. In the more general case of a (non-discrete) 
locally compact group, the modular function can be non-trivial.
Then the property $h\delta=h$ is related with the fact that
compact groups are unimodular. We refer to [L-VD2] with non-trivial examples for an illustration of the above results.
\snl
In most of the paper, we will deal with a single group-like
projection $h$. Then we can {\it normalize} the left integral
$\varphi$ and the right integral so that $\varphi(h)=1$ and
$\psi(h)=1$. This is possible because we have that $\varphi(h)>0$
and $\psi(h)>0$ in general. Because the scaling constant $\nu$
equals $1$, this joint normalization of the left and the right
integral will give that $\psi=\varphi \circ S$. So, in most of the
paper, when we only deal with one specific group-like projection
$h$, we will assume that $\varphi(h)=\psi(h)=1$.
\snl
On a few occasions, we will treat several group-like projections
together. This will be done e.g.\ furhter in this section. In that
case, we will not assume the normalizations above. We will clearly
mention this fact when it appears.

\nl
\it The Fourier transform of a group-like projection \rm

\nl
Now we want to study the Fourier transform of a group-like
projection. The remarkable property is that it is again a
group-like projection.
\snl
First we recall the framework, as it is reviewed in the
appendix, needed to consider the Fourier transform.
\snl
So, as before,
$(A,\co)$ is a multiplier Hopf \st algebra with positive
integrals, i.e.\ an algebraic quantum group. We use $B$ for the
dual $\hat A$ and we take for the comultiplication $\co$ on $B$
the opposite comultiplication $\hat\co^{\text{op}}$ of the dual
comultiplication $\hat\co$ on $\hat A$. We also use
$\langle\,\cdot\,,\,\cdot\,\rangle$ for the pairing between $A$
and $B$. Observe that one of the consequences of this convention
is that $\langle S(a),b\rangle=\langle a,S^{-1}(b)\rangle$ when
$a\in A$ and $b\in B$. As we explain in the appendix, this convention is very common in the theory of
locally compact quantum groups.
\snl
For the Fourier transform, we make the choice
$F(a)=\varphi(\,\cdot\, a)$ where $\varphi$ is a left integral on $A$. We
know that then the inverse transform is given by
$a=\varphi(S^{-1}(\,\cdot\,)b)$ when $b=\varphi(\,\cdot\, a)$, provided
we have the correct relative normalization of the left integrals
on $A$ and on $B$ (Proposition 3.4 in [VD5]). Again, see the appendix for more details about this Fourier transform.
\snl
We have the following result.

\inspr{1.8} Proposition \rm
Let $h$ be a group-like projection in the algebraic quantum group
$(A,\co)$. Put $k=F(h)=\varphi(\,\cdot\,h)$. Then $k$ is a
group-like projection in the dual $B$.

\snl\bf Proof: \rm
Recall that we assume $\varphi(h)=1$. 
\snl
We first show that $k^2=k$. We will again make use of the Sweedler notation for convenience. For any $a\in A$ we have
$$\align
      \langle a, k^2 \rangle
           &= \langle a_{(1)}, k \rangle \langle a_{(2)}, k \rangle 
           =   \varphi(a_{(1)}h) \varphi (a_{(2)}h)  \\
           &=  (\varphi\ot\varphi)(\co(a)(h\ot h))  
           =  (\varphi\ot\varphi)(\co(a)\co(h)(h\ot 1))  \\
           &=  \varphi(ah) \varphi(h) =  \varphi(ah) \\
           &= \langle a, k \rangle.
\endalign$$
Next we show that $k^*=k$. 
By the definition of the adjoint on the
dual, we have $\langle a , k^* \rangle = \langle S(a)^* , k
\rangle ^- =  \varphi(S(a)^*h)^- $ for any $a \in A$. As $\varphi$
is positive, hence self-adjoint, we get $\varphi (S(a)^*h)^- =
\varphi (h S(a))$. Because $S(h)=h$ we have $\varphi
(hS(a))=\varphi(S(ah))$. Using the property
$\varphi(S(x))=\varphi(x\delta)$ for all $x\in A$, we get
$\varphi(S(ah))=\varphi(ah\delta)$. In Proposition 1.7 we found
$h\delta =h$. Combining all these equations, we get $\langle a,
k^* \rangle = \varphi(ah) = \langle a , k \rangle $ for all $a$ so
that $k^* = k$. 

\snl 
Finally, we show that $\co(k)( 1 \ot k) = k \ot k$. In the next series of equalities, we use that the product in $B$ is dual to the coproduct in $A$, as well as many of the results about $h$ from
Propositions 1.6 and 1.7. We will need $\co(h)(h \ot 1) = h \ot h$, $S(h)=h$, $h\delta=h$ and $\sigma(h)=h$. For any pair $a,a' \in A$ we have
$$\align
    \langle a \ot a' , \co(k)(1\ot k) \rangle
  &= \langle a \ot a'_{(1)} , \co(k) \rangle \langle  a'_{(2)},k\rangle 
           = \langle  a'_{(1)}a, k \rangle \langle  a'_{(2)} , k \rangle  \\
  &= (\varphi \ot \varphi) (\co(a')(a \ot 1)(h\ot h)) \\
  &= (\varphi \ot \varphi) ((h \ot h)\co (a')(a \ot 1)) \\
  &= (\varphi \ot \varphi) ((h \ot 1) \co (ha') (a \ot 1)) \\
  &= \varphi(ha') \varphi(ha) 
            = \varphi(a'h) \varphi(ah) \\
  &= \langle a \ot a' , k \ot k \rangle.
\endalign$$
\einspr

Of course, the formulas that we have proven for $h$ in the
Propositions 1.6 and 1.7 will also be valid for $k$. And in fact, it is
interesting to see how these various properties of $k$ relate with
the various properties of $h$. In [L-VD2], where we consider duality, we will see some of these results.
We will e.g.\ show that $h$ and $k$ commute in the Heisenberg algebra.

\nl
At this point, it is interesting to have a look at {\it the case of two group like projections}, one smaller than the other.
\snl
As before, let $(A,\co)$ be an algebraic quantum group and let $h$
and $h'$ be two group-like projections in $A$. Now, we will not
assume any specific normalization of the integrals. We consider
the case $h\leq h'$. In this algebraic setting, the inequality
$h\leq h'$ means that $hh'=h'h=h$. Observe that one equation
$hh'=h$ implies the other $h'h=h$ by taking adjoints.
\snl
Using that $h'$ is group-like and $h'h=h$, we get immediately
$\Delta(h')(1\ot h)=h'\ot h$. Similarly $\Delta(h')(h\ot 1)=h\ot h'$. Then, we also
get the following property.

\inspr{1.9} Proposition  \rm
As before, let $h$ and $h'$ be group-like projections
in an algebraic quantum group. Consider the (properly scaled)
Fourier transforms $k$ and $k'$ of $h$ and $h'$ respectively. If
$h\leq h'$ then $k'\leq k$.

\snl\bf Proof: \rm
Denote $c=\varphi(h)$ and $c'=\varphi(h')$. Then
$k=c^{-1}\varphi(\,\cdot\, h)$ and $k'={c'}^{-1}\varphi(\,\cdot\,
h')$. For any $a\in A$ we get
$$\align cc'\langle a, kk' \rangle
   &=\varphi(a_{(1)}h)\varphi(a_{(2)}h') \\
   &=\varphi(a_{(1)}h'_{(2)}S^{-1}(h'_{(1)})h)\varphi(a_{(2)}h'_{(3)}).
\endalign$$
Because $S(h)=h$ we get
$$\align S^{-1}(h'_{(1)})h \ot h'_{(2)}
     &=(S^{-1}\ot\iota)((h\ot 1)\co(h')) \\
     &=(S^{-1}\ot\iota)(h\ot h')=h\ot h'.
\endalign$$
Hence
$$\align cc'\langle a, kk' \rangle 
    &=\varphi(a_{(1)}h'_{(1)}h)\varphi(a_{(2)}h'_{(2)})\\
    &= \varphi(h)\varphi(ah')=cc' \langle a, k' \rangle.
\endalign$$
So $kk'=k'$ and therefore $k'\leq k$.
\einspr

We see that increasing projections $h\leq h'$ yield decreasing Fourier transforms $k'\leq k$.
\nl
\it Examples and some special cases \rm
\nl
It is not hard to illustrate the results of Proposition 1.8 in the
case of our motivating examples. Take the example of a discrete
group (Example 1.2.i and 1.2.ii). These two cases are dual to each
other. Using the pairing $\langle f , \lambda_p \rangle = f(p)$,
we see that $k$, defined as $\frac{1}{n} \sum_{p \in k} \lambda_p$ (with $n=\# K$), and
$h$, defined as $\chi_K$, satisfy
$$\langle f , k \rangle = \frac{1}{n}\sum_{p \in K} f(p) = 
\varphi(fh).$$
We use $\varphi(f)=\frac{1}{n}\sum_{p\in G} f(p)$ (taking care of the correct normalization of $\varphi$).

\nl
The last thing we will do in this section is to consider two
special cases: i) $h$ is central in $A$ and ii) $k$ is central in
$B$. In the case of the example i) in 1.2, of course $A$ is
abelian and $h$ is automatically central. When the group is not
abelian, $B$ is not abelian. It turns out that $k$ being central
is related with the fact that $K$ is a normal subgroup. The reader
should have this example in mind for a better understanding of some aspects of the rest of the paper.

\snl
We first formulate (and prove) the main result in this connection. Then we give some comments. 

\iinspr{1.10} Proposition \rm
With the notations as before, the following are equivalent:
\item{}\ i) $h$ is central in $A$,
\item{}\ ii) $hA=Ah$,
\item{}\ iii) $\co^{\text{op}}(k) = (\iota \ot S^2)(\co(k)(1\ot \delta^{-1}))$,
\item{}\ iv) $\{\langle a, k_{(1)}\rangle  k_{(2)} \mid a \in A \}
             =\{\langle a, k_{(2)}\rangle  k_{(1)} \mid a \in A \}$.

\snl\bf  Proof: \rm 
Let us first show that i) and ii) are
equivalent. It is clear that i) implies ii) and we only have to
show the converse. So, assume that $hA=Ah$. Take any $a \in A$.
Then there exists an $a'\in A$ so that $ha=a'h$. Then
$$ha=a'h=a'hh=hah.$$
So $ha=hah$ for all $a\in A$. Taking adjoints, we get that also
$ah=hah$ for all $a$. Hence $ha=ah$ for all $a \in A$. Therefore also
ii) implies i).
\snl
Now we will prove that i) and iii) are equivalent. Take a pair of elements $a,a' \in
A$. Using that $\varphi$ is $\sigma$-invariant, we get
$$\varphi(a'\sigma(a)h)=\varphi(\sigma^{-1}(a')ah)=
    \langle\sigma^{-1}(a')a , k \rangle=
    \langle a \ot \sigma^{-1}(a'), \co(k) \rangle.$$
Now we know from the appendix that
$$\langle \sigma^{-1}(a),b \rangle = \langle a, S^2(b)\delta^{-1}
    \rangle$$
for all $a\in A$ and $b\in B$. If we use this formula in the previous one, we get
$$\varphi(a'\sigma(a)h)=
   \langle a \ot a' , (\iota \ot S^2) \co(k)(1 \ot
   \delta^{-1})\rangle.$$
On the other hand, we also have
$$\langle a \ot a' , \co^{\text{op}}(k) \rangle = \langle aa', k
    \rangle = \varphi(aa'h)=\varphi(a'h\sigma(a)).$$
So we see that condition iii) is equivalent with the equality
$$\varphi(a'h\sigma(a))=\varphi(a'\sigma(a)h)$$
for all $a,a' \in A$. This is of course true when $h$ is central. So
i) implies iii). Conversely, if iii) holds, then this equaltiy is true
for all $a,a' \in A$. By the faithfulness of $\varphi$, it then
follows that $h\sigma(a)=\sigma(a)h$ for all $a\in A$. Then $h$ is
central. This proves the equivalence of i) and iii)
\snl
Finally, let us show that ii) and iv) are equivalent. From the
definition of $k$ (see also the previous calculations), we get
$$\{\langle a , k_{(1)}\rangle k_{(2)} \mid a \in A\}=
  \{\varphi(\,\cdot\, ah) \mid a \in A \}.$$
Similarly
$$\align
   \{\langle a , k_{(2)}\rangle k_{(1)} \mid a \in A\}
     &=\{\varphi(a\,\cdot\, h) \mid a \in A \}. \\
     &=\{\varphi(\,\cdot\, h\sigma(a)) \mid a \in A \} \\
     &=\{\varphi(\,\cdot\, ha) \mid a \in A \}.
\endalign$$
So, condition iv) will hold if and only if $Ah=hA$, that is if ii)
holds.
\einspr

Condition iv) simply means that 'the right leg' of $\co(k)$ is
the same as its 'left leg' (see the appendix for some more details about how to define these 'legs').
Condition iii) looks strange, but, as we saw in the proof, it is the dual version of condition i). So, the equivalence of iii) and iv) is as the equivalence between i) and ii).

\snl
By duality, we also have this result when we replace $h$ by $k$. So, we get e.g.\
that $k$ is central if and only if the two legs of $\co(h)$ coincide.

\snl
In the case of a group, we have $S^2=\iota$ and also $\delta=1$ in
$B$. Then $h$ is central and $\co(k)$ is symmetric. In this case,
$k$ will be central if the left leg of $\co(h)$ coincides with its
right leg. Now, recall that $h$ is the characteristic function of
the subgroup $K$. These two legs of $\co(h)$ are nothing else but
the functions that are constant on right, respectively left cosets.
These are the same if and only if $K$ is a normal subgroup. 

\snl
We will come back to these various cases later, especially, in the beginning
of Section 2 and of Section 3.

\nl\nl

\bf 2. Compact quantum hypergroups \rm
\nl
As before, let $(A,\co)$ be a multiplier Hopf \st algebra and
assume that $h$ is a group-like projection (cf.\ Definition 1.1).
So $h$ is a non-zero element in $A$ satisfying $h^2=h=h^*$ and
$\co(h)(1 \ot h)=h\ot h$. We have seen that, in the case of an
algebraic quantum group, automatically also $\co(h)(h \ot 1)=
h\ot h$ and that $\varepsilon(h)=1$ and $S(h)=h$ (cf.\ Proposition 1.6). 
In any case, let us assume that also $\co(h)(h \ot 1)= h \ot h$. We know that 
then the two other formulas, namely $\varepsilon(h)=1$ and $S(h)=h$, are
also true.
\snl
In this section, we will first consider the case where $h$ is a central 
projection. Then, we will see that $Ah$ is a compact quantum
group. In this paper, by a compact quantum group, we mean a Hopf $^*$-algebra 
with a positive integral (i.e.\ an algebraic quantum group of compact type in
the sense of [VD3]). This is essentially the algebraic ingredient of a compact
quantum group as defined in the C$^*$-algebraic sense by Woronowicz (see [W]). 
In the general case, when $h$ is not a central projection, we consider
$hAh$  and 
we encounter a compact quantum hypergroup (in the sense of [De-VD1]). We use  
the basic examples to illustrate the results. They turn out to give the classical Hecke algebras.
\nl
\it The case of a central group-like projection \rm
\nl
The following result is fairly easy to obtain.

\inspr{2.1} Proposition \rm Assume that $h$ is central in $A$ and
denote $A_0=Ah$. Then $A_0$ is a Hopf \st algebra when we define
the coproduct $\co_0$ on $A_0$ by $\co_0(a)=\co(a)(h \ot h)$.

\snl \bf Proof: \rm
First observe that $\co_0(A_0) \subseteq A_0 \ot A_0$. Indeed,
$$\co(a)(h \ot h)=(\co(a)(1 \ot h))(h \ot h)$$
and $\co(a)(1 \ot h) \in A \ot A$. Next remark that $\co_0$ is a
\st homomorphism from $A_0$ to $A_0 \ot A_0$, precisely because
$h$ is central in $A$. Of course, $h$ is the identity in $A_0$ and
$\co_0(h)= \co(h)(h\ot h)=h\ot h$ so that $\co_0$ is unital.
\snl
To prove coassociativity of $\co_0$, let $a \in A_0$ and write
$$\align
    (\co_0 \ot \iota)\co_0(a)
       &= ((\co \ot \iota)\co_0(a))(h \ot h \ot 1) \\
       &= (\co \ot \iota)(\co(a)(h \ot h))(h \ot h \ot 1) \\
       &= ((\co \ot \iota)\co(a))(\co(h)(h \ot h) \ot h) \\
       &= ((\co \ot \iota)\co(a))(h \ot h \ot h).
\endalign$$
Similarly
$$(\iota \ot \co_0)\co_0(a)= ((\iota\ot\co)\co(a))(h\ot h\ot h)$$
and the coassociativity of $\co_0$ follows from the coassociativity of
$\co$.
\snl
We also see that 
$$\align
    \co_0(A_0)(1\ot A_0)
       &=(\co(A)(1\ot A))(h\ot h) \\
       &=(A\ot A)(h\ot h) \\
       &=A_0 \ot A_0
\endalign$$
and similarly
$$(A_0 \ot 1)\co_0(A_0)=A_0 \ot A_0.$$ Finally, it is easily verified that the linear maps $a\ot b\mapsto \co_0(a)(1\ot b)$
and $a\ot b\mapsto (a\ot 1)\co_0(b)$ are injective on $A_0\ot A_0$.
It follows that $A_0$ is a Hopf \st algebra (see e.g.\ [VD1]).
\einspr

It is not so hard to verify that the counit $\varepsilon_0$ for
$\co_0$ is simply the restriction of $\varepsilon$ to $A_0$
because $\varepsilon(h)=1$. Similarly, the antipode $S_0$ is given
by the restriction of $S$ to $A_0$ as $S(h)=h$.

\nl
If $(A,\co)$ is a multiplier Hopf $^*$-algebra with integrals, we get the
following.

\inspr{2.2} Theorem \rm
Let $(A,\co)$ be an algebraic quantum group. Assume that $h$ is a
central group-like projection. Then $(A_0,\co_0)$ is a compact
quantum group when, as in the previous proposition, $A_0=Ah$ and
$\co_0$ is defined by $\co_0(a)=\co(a)(h\ot h)$. The Haar functional
on $A_0$ is given by the restriction of the left integral on $A$
to $A_0$.
\einspr

Observe that the property $h\delta=h$, proven in Proposition 1.7,
gives that the left and the right integrals on $A$ have the same
restriction to $A_0$ (as expected).
\snl
The typical example to illustrate the above result is coming from
a compact open subgroup $K$ of a locally compact group $G$. In
this case, as we saw in Example 1.3.i, $A$ is the \st algebra of
'polynomial functions' on $G$ and $h$ is the characteristic
function of $K$. As $A$ is abelian, $h$ will be central and
$(A_0,\co_0)$ is nothing else but the Hopf \st algebra of
polynomial functions on the compact group $K$.
\nl
\it The general case \rm
\nl
In the general case, when $h$ is no
longer assumed to be central, we have a far more complicated
situation.
\nl
To begin with, $Ah$ is no longer a \st algebra. If it is a \st
algebra, then $Ah=hA$ and we have seen in Proposition 1.10 that
then $h$ is central.
\snl
There are two ways to get a \st algebra. On can either take $hAh$ or
one can consider $AhA$. But only on the first one, it seems to be
possible to define some kind of a comultiplication. Now, we must
restrict on both sides by $h\ot h$ as we show in the following proposition.

\inspr{2.3} Proposition \rm
Let $A_0=hAh$. Define $\co_0$ on $A_0$ by
$$\co_0(a)=(h\ot h)\co(a)(h\ot h).$$
Then $\co_0$ is a unital, positive map from $A_0$ to
$A_0 \ot A_0$ and it is coassociative.

\snl\bf Proof: \rm
This is all more or less obvious. The proof of the coassociativity is
essentially the same as in Proposition 2.1. The positivity can be seen
as follows:
$$\align
     \co_0(a^*a)
        &=(h\ot h)\co(a)^*\co(a)(h\ot h) \\
        &=(\co(a)(h\ot h))^*(\co(a)(h\ot h))
\endalign$$
for any $a\in A_0$.
\einspr

One can show that $\co_0$ is not only positive, but also completely
positive. Recall that {\it complete positivity} for
$\co_0$ means that $\co_0 \ot \iota$ defined on $A_0\ot M_n(\Bbb
C)$ is still positive for all $n$ where $M_n(\Bbb C)$ is the \st
algebra of $n$ by $n$ complex matrices with the usual \st
operation.  Because it is not essential here, we will not give
details. 
\snl
The counit $\varepsilon$ can be restricted to $A_0$ and will give
a unital \st homomorphism $\varepsilon_0$ from $A_0$ to $\C$
satisfying the usual conditions:
$$(\varepsilon_0 \ot
       \iota)\co_0(a)=a=(\iota\ot\varepsilon_0)\co_0(a)$$
for all $a\in A_0$.
\snl
The antipode $S$ can also be restricted to $A_0$ because
$S(h)=h$. This restriction $S_0$ is a anti-homomorphism of $A_0$ and
it satisfies $S_0(S_0(a)^*)^*=a$ and
$$\co_0(S_0(a))=(S_0\ot S_0)\co^{\text{op}}_0(a)$$
for all $a \in A_0$. Recall that $\co^{\text{op}}_0$ is the
opposite comultiplication, obtained by composing $\co_0$ with the
flip. Unfortunately, these conditions are too weak to fully
characterize the antipode. If e.g.\ $A_0$ is abelian and
coabelian, the identity map will also satisfy these formulas.

\nl
Let us now proceed in the case of an {\it algebraic quantum
group}. We get the following results about the integrals. We
denote by $\varphi_0$ the restriction of the left integral
$\varphi$ to the algebra $A_0$.

\inspr{2.4} Proposition \rm
Assume that $(A,\co)$ is an algebraic quantum group and that $h$ is a
group-like projection in $A$. Let $A_0$ and $\co_0$ be as in
Proposition 2.3. The left integral $\varphi$ and the right integral $\psi$
coincide on $A_0$. The restriction $\varphi_0$ of 
$\varphi$ to $A_0$ is a positive linear functional on $A_0$ satisfying
$$\align
    (\iota\ot\varphi_0)\co_0(a) &= \varphi_0(a)h \\
    (\varphi_0\ot\iota)\co_0(a) &= \varphi_0(a)h
\endalign$$
for all $a\in A_0$.

\snl\bf Proof: \rm
This is essentially straightforward. For the second equatity, we
use that the left integral $\varphi$ and the right integral $\psi$
coincide on $A_0$. This is so because $h\delta=h$.
\einspr

Remark that $h$ is the identity in $A_0$ and so these are the usual
formulas to express left and right invariance of the functional
$\varphi_0$. Observe also that $\varphi_0(h)=\varphi(h)$ and that this
is non-zero. It follows that left and right invariant functionals are
unique (and equal).
\snl
In the case of an algebraic quantum group, the antipode on the pair $(A_0,\co_0)$ can be
characterized using the relation with the invariant integral in the
following way.

\inspr{2.5} Proposition \rm
We assume again that $h$ is a group-like projection in the algebraic
quantum group $(A,\co)$. We will use the notations of before. So
again, $S_0$ is the restriction of the antipode to the \st subalgebra
$A_0$. Also for this restriction we have
$$S_0((\iota\ot\varphi_0)(\co_0(a)(1\ot b)))
    =(\iota\ot\varphi_0)((1\ot a)\co_0(b))$$
for all $a,b\in A_0$.

\snl\bf Proof: \rm
For $a,b\in A_0$ we have
$$(\iota\ot\varphi_0)(\co_0(a)(1\ot b))
       = (\iota\ot\varphi)((h\ot h)\co(a)(h\ot h)(1\ot b)).$$
This is $hxh$ where $x=(\iota\ot\varphi)((1\ot h)\co(a)(1\ot hb))$. We get 
$$x = (\iota\ot\varphi)(\co(a)(1\ot hbh)) 
      = (\iota\ot\varphi)(\co(a)(1\ot b))$$
because $\sigma(h)=h$ and $hbh=b$. We know that
$$S(x)=(\iota\ot\varphi)((1\ot a)\co(b))$$
(see e.g.\ Proposition A.4 in the appendix). A similar calculation will give that
$$hS(x)h=(\iota\ot\varphi_0)((1\ot a)\co_0(b)).$$
So, because $S(h)=h$, we get the desired formula.
\einspr

We know that the formula in the previous proposition characterizes the
antipode in the case of an algebraic quantum group. However, we need
to know that all elements in $A_0$ can be written as (a linear
combination of) elements of the form
$$(\iota\ot\varphi_0)(\co_0(a)(1\ot b))$$
where $a,b\in A_0$. This is true for $(A,\co)$. Indeed, the linear span of
$$\{(\iota\ot\varphi) (\co(a)(1\ot b)) \mid a,b \in A \}$$
is $A$ simply because $\co(A)(1\ot A)=A\ot A$. This stronger condition will not remain valid for the
pair $(A_0,\co_0)$ however. In general, it will not be true that
$\co_0(A_0)(1\ot A_0)=A_0 \ot A_0$. See the examples below.
\snl
On the other hand, we do have the weaker condition. This is the content of 
the next proposition.

\inspr{2.6} Proposition \rm
We have
$$A_0=\text{span}\{(\iota\ot\varphi_0)((\co_0(a)(1\ot b)) \mid a,b\in
A_0\}.$$
\snl\bf Proof\rm:
If the result is not true, then there is a non-zero linear functional 
$\omega$ on $A_0$ that kills all the elements of the right hand side. By the 
fact that $\varphi_0$ remains faithful on the subalgebra $A_0$, it follows that
$(\omega\ot\iota)\Delta_0(a)=0$ for all $a\in A_0$. Now apply the counit $\varepsilon_0$
to get $\omega(a)=0$ for all $a$. This is a contradiction. 
\einspr

If we combine all these results, we get the following.

\inspr{2.7} Theorem \rm
The pair $(A_0,\Delta_0)$ is an algebraic quantum hypergroup of compact type
(in the sense of [De-VD1]).
\einspr

Recall that it is called of compact type because the algebra $A_0$ has an 
identity. We get an algebraic version of a compact
quantum hypergroup as it was introduced by Chapovsky and Vainerman, see e.g. [C-V] (and some earlier references therein). 
\snl
In the next section, we will encounter an algebraic
quantum hypergroup of discrete type (see Theorem 3.17).
\nl
In the previous section, we have looked at the situation of a pair of group-like projections
$h,h'\in A$ such that $h\leq h'$. For the results in this section, this situation is
not really very interesting because the corresponding algebra $hAh$ is
simply a subalgebra of  $h'Ah'$ (because
$h=h'h=hh'$. Remark however that the coproduct on the smaller one is
not just the restriction of 
the coproduct on the bigger one. One has to cut down further with $h\otimes h$ on each side. 
\nl
\it Examples \rm
\nl
Let us now look at our examples and see what we get for the pair $(A_0,\Delta_0)$ 
in these simple special cases. 

\inspr{2.8} Examples \rm
i) First take the case of a (discrete) group $G$ with a finite subgroup $K$. Let
$A$ be the group  algebra $\Bbb C G$ and put $h=\frac{1}{n} \sum_{p\in
K} \lambda_p$ as in Example 1.2.ii.
\snl
When the group $K$ is a normal subgroup, we have, for any $q\in G$,
$$\lambda_q h=\frac{1}{n}\sum_{p\in K} \lambda_{qp}
     =\frac{1}{n} \sum_{p\in K}
     \lambda_{qpq^{-1}}\lambda_q=h\lambda_q$$
and $h$ is central. Here $Ah$ is the group algebra of the quotient
group $G/K$.
\snl
ii) Again, let $K$ be a finite subgroup of $G$, not necessarily
normal. Let $A$ and $h$ be as in i).
For any $q\in G$ we have
$$h\lambda_q h = \frac{1}{n^2} \sum_{p,p'\in K} \lambda_{pqp'}.$$
Let us denote this element by $\pi(q)$ and observe that it only
depends on the double coset $KqK$. For the product of such elements,
we have
$$\align
    \pi(q)\pi(q')
       &= \frac{1}{n^4}\sum_{p,p',r,r'} \lambda_{pqp'rq'r'} \\
       &= \frac{1}{n^3}\sum_{p,p',r} \lambda_{pqrq'p'} \\
       &= \frac{1}{n} \sum_{r\in K} \pi(qrq').
\endalign$$
For the coproduct $\co_0$, we get
$$\align
    \co_0(\pi(q))
       &= (1 \ot h)\co(\pi(q))(1\ot h) \\
       &= \frac{1}{n^4}\sum_{p,p',r,r'} \lambda_{pqp'} \ot
                       \lambda_{rpqp'r'} \\
       &= \frac{1}{n^4}\sum_{p,p',r,r'} \lambda_{pqp'} \ot
                       \lambda_{rqr'} \\
       &= \pi(q) \ot \pi(q).
\endalign$$
The counit $\varepsilon_0$ is given by the restriction of
$\varepsilon$ and so $\varepsilon(\pi(q))=1$ for all $q$. The antipode
$S_0$ is given by the restriction of $S$ and so
$S_0(\pi(q))=\pi(q^{-1})$. Observe that in general,
$$\pi(q^{-1})\pi(q)\neq 1$$
so that the formula $m(S_0\ot \iota)\co_0(a)=\varepsilon_0(a)1$,
true in a Hopf algebra, will not hold here. Similarly it will not be true that
$\co_0(A_0)(1\ot A_0)=A_0 \ot A_0$.
\snl
For the integral, we have $\varphi(\lambda_p) = 0$ except for
$p=e$, the identity in the group. So, when we restrict, we will
get $\varphi_0(\pi(q))=0$ when $q\notin K$. Remark that $\pi(q)=h$
if and only if $q\in K$. In that case of course, $\varphi_0(\pi(q))\neq 0$.
\snl
Finally, let us look at the formula in Proposition 2.5. And first, let
us calculate
$$(\iota\ot\varphi_0)(\co_0(\pi(q))(1\ot \pi(q'))$$
for any two elements $q,q'\in K$.
We get
$$\pi(q)\varphi_0(\pi(q)\pi(q'))
         =\pi(q) \frac{1}{n} \sum_{r\in K} \varphi_0(\pi(qrq')).$$
This is $0$ except when $\pi(q')=\pi(q^{-1})$ and in that case, it
is $\pi(q)$, provided we use the normalization $\varphi(h)=1$.
\snl
Similarly,
$$(\iota\ot\varphi_0)(1\ot \pi(q))\co_0(\pi(q'))=0$$
except again when $\pi(q')=\pi(q^{-1})$, in which case the
expression is $\pi(q')$. We see that these elements indeed span
$A_0$ and that the antipode $S_0$ satisfies the equality of Proposition 2.5.
\einspr

The case of Example 1.2.i. is too trivial to mention. 
\snl
Let us now also consider the example with a  compact open subgroup of
a locally compact group. Again Example 1.3.i. is not very interesting and also here, we only look at the second case in Example 1.3.

\inspr{2.9} Example \rm
Let $G$ be a locally compact group with a compact open subgroup $K$. Consider the reduced C$^*$-algebra $C_r^*(G)$ and its dense multiplier Hopf $^*$-algebra $A$ of polynomial functions (as in Example 1.3.ii). Take for $h$ the characteristic function of $K$, as sitting in $A$. We have seen that $h$ is a group-like projection. 
\snl
If $K$ is a normal subgroup, the element $h$ is central and we end up by considering the group algebra of the (disrete) quotient group $G/K$.
\snl
In general, when $K$ is not assumed to be a normal subgroup, we must consider the algebra $A_0$, defined as $hAh$. This 
consists of elements of the form $\int_G f(p) \lambda_p dp$ where $f$ is a continuous function with compact support on $G$, constant on double cosets $K\backslash G/K$. The coproduct $\Delta_0$ on $hAh$ is given by the formula
$$\Delta_0 (\int_G f(p) \lambda_p dp)= \frac1{n^2}\int_G  \int_K\int_K f(p)\,\lambda_p \ot \lambda_{kpk'}\, dk\,dk'\,dp$$
when $f\in hAh$ where as before, $n$ is the measure of $K$. The counit and antipode on $A_0$ are the restrictions of the counit and antipode on $A$ and the same is true for the integral $\varphi_0$ on $A_0$.
\snl
In fact, because the space of double cosets will be a discrete space, we get a compact quantum hypergroup of the same type as in the previous example. Indeed, it is still possible to consider the elements $\pi(q)$ for any $q\in G$. Now they are defined, similarly as in Example 2.8.ii by 
$$\pi(q)=h\lambda_q h= \frac1{n^2}\int_K\int_K \, \lambda_{kqk'}\, dk\,dk'$$
and $A_0$ is again the linear span of these element $\pi(q)$ with $q\in G$. 
As before, we will get $\Delta_0(\pi(q))=\pi(q) \ot \pi(q)$ for all $q\in G$. Also $\varepsilon_0(\pi(q))=1$ and $S_0(\pi(q))=\pi(q^{-1})$ for all $q$. Finally, also $\varphi_0(\pi(q))=0$ except for $q\in K$. 
\einspr

We finish this section with some remarks.

\iinspr{2.10} Remarks \rm
i) As we saw in Proposition 2.3, the situation is quite nice, and
in fact not very surprising, when $h$ is a {\it central}
group-like projection. When the projection is not central, we have
considered the \st algebra $hAh$ as this algebra still carries a
reasonable structure (it is a compact quantum hypergroup).
\snl
ii) In the Example 2.8.ii, where we have a (discrete) group $G$ with a finite subgroup $K$, the algebra we get is nothing else but the classical Hecke algebra $H(G,K)$ (see e.g.\ Section 2.10 in [G-H-J]). Also the Example 2.9 gives an algebra of this type.
\snl
iii) In general however, this algebra may be rather small and contain too little
information. Consider e.g.\ the case of $G=GL(2,\Bbb Q)$ with the
subgroup $K=SL(2,\Bbb Z)$. Consider the group algebra (as in
Example 2.9) so that $A_0$ is the convolution
algebra with functions constant on double cosets. In this example,
this algebra turns out to be abelian and it will not contain much
information about the pair $(G,K)$. See e.g.\ [Kr].
\snl
iv) Finally, let us make a short remark about terminology. In the Examples 1.2.i and 1.3.i, the group-like projection $h$ is central as the algebras are abelian. We have that $A$ is an algebra of functions on the group $G$ and $Ah$ is an algebra of functions on the subgroup $K$. From this point of view, in the general case, we can think of $Ah$ (if $h$ is central) and $hAh$ (in general) as a compact quantum (hyper) {\it sub}group. However, if we consider  Example 2.8 (cf.\ Example 1.2.ii) and Example 2.9 (cf.\ Example 1.3.ii), we see that $A$ is a group algebra and, in the case of a normal subgroup, that $Ah$ is a group algebra of the quotient group. From this point of view, in the general case, we can think of $Ah$ (if $h$ is central) and $hAh$ (in general) as a compact quantum (hyper) {\it quotient} group. 
\einspr

In [L-VD2] we plan to give some other, perhaps more interesting examples of compact quantum hypergroups obtained in this way.
\nl\nl

\bf 3. Discrete quantum hypergroups \rm
\nl
In this section, we will study another \st algebra associated with a
group-like projection. The case is dual to the case considered in the
previous section as we show in the second paper on this subject ([L-VD2]. We will consider the dual objects to $Ah$, $hA$ and $hAh$.
\snl
In the previous section, we got a compact quantum group if $h$ is central. In this section, we first look at
the case where the two legs of $\Delta(h)$ coincide. This will yield a 
discrete quantum group. Recall that in this paper, a discrete quantum group is 
an algebraic quantum group with co-integrals (i.e.\ an algebraic quantum
group of discrete type as in [VD3], see also [VD2]). In the general case, we will encounter
discrete quantum hypergroups (in the sense of [De-VD1]). Again, we also consider
the basic examples. In this case, the underlying algebras consist of functions
that are constant on cosets (whereas in the previous section, we rather looked
at the functions living on the subgroup). 
\snl
We begin with considering the algebra $C$. The reader should have in mind that this is dual to $Ah$ and in fact the (inverse) Fourier transform of the space $Bk$ (see again [L-VD2] for more details about this point of view).
\nl
\it The left invariant \st subalgebra $C$ \rm
\nl
Again, we start with the general situation and see what is
possible. So, to begin with, let $(A,\co)$ be any multiplier Hopf
\st algebra and $h$ a group-like projection in  $A$. We will assume
that $(A,\co)$ is regular.
\snl
In what follows, we will use the definitions and results of the
appendix concerning the legs of the coproduct (cf.\ Definition A.11 and following).

\inspr{3.1} Notation \rm
Let $C$ denote the right leg of $\co(h)$. If we work with several
group-like projections (see further in this section), we will write $C_h$. But we
can drop the index because most of the time, we only work with a single
group-like projection.
\einspr

By the self-adjointness of $\co(h)$, we certainly have that $C$ is a \st
subspace of $A$. But we can prove more. We have the
following properties of elements in $C$.

\inspr{3.2} Proposition \rm
For any $a\in C$ we have
\smallskip
\item{} i) $ah=\varepsilon(a)h$,
\item{} ii) $\co(a)(1\ot h)=a\ot h$,
\item{} iii) $(S(a)\ot 1)\co(h)=(1\ot a)\co(h)$.
\snl\bf Proof: \rm
Start from the equation $\co(h)(1\ot h)=h\ot h$. If we apply a linear functional
$\omega$ on the first leg of this equation, we get $ah=\omega(h)h$
where $a=(\omega\ot\iota)\Delta(h)$. 
Clearly $\omega(h)=\varepsilon(a)$ and because every element in $C$ is
of this form, we have proven
 i).
\snl
To prove ii), start again with $\co(h)(1 \ot h)=h \ot h$, apply
$\co\ot\iota$, use coassociativity and multiply with $a\ot 1\ot 1$ to
get
$$((\iota\ot\co)(\co(h)(a\ot 1))(1\ot 1\ot h) = (\co(h)(a\ot 1))\ot
h.$$
Now take a linear functional $\omega$ on $A$ and apply
$\omega\ot\iota\ot\iota$ to get $\co(x)(1\ot h) = x\ot h$
where $x=(\omega\ot\iota)(\co(h)(a\ot 1))$. Again, we get all elements $x$ in
$C$  and so we proved ii).
\snl
To prove iii), we start with $\co(a)(1\ot h )=a\ot h$, where $a\in C$,
and we apply $\iota\ot\co$. We will also use the Sweedler notation and write
$$(\Delta\ot\iota)\Delta(a)= a_{(1)} \ot a_{(2)} \ot a_{(3)}.$$ 
Then
$$ a_{(1)} \ot (a_{(2)} \ot a_{(3)})\co(h) = a \ot \co(h).$$
Now apply $S$ on the first leg of this equation and multiply to get
$$ (S(a_{(1)}) a_{(2)} \ot a_{(3)})\co(h) = (S(a) \ot 1) \co(h).$$
The left hand side is $(1\ot a)\co(h)$ and so we get the result.
\einspr

Observe that we have proven condition iii) from condition ii)
without the need for any other property of $h$. This fact will be
used later, in the proof of Proposition 3.4.

\inspr{3.3} Proposition \rm
The space $C$ is a \st algebra and $\co(C)(A\ot 1)\subseteq A\ot
C$. We have $h\in C$.
\snl\bf Proof: \rm
We have shown already that $C$ is \st invariant. To show that it is a
subalgebra, take any $a\in C$ and $a'\in A$. Then, from iii) in the
previous proposition, we get
$$(1\ot a)\co(h)(a'\ot 1) = (S(a) \ot 1)\co(h)(a' \ot 1) \subseteq
A\ot C.$$
Take an element $\omega\in A'$ and apply $\omega\ot\iota$ to this
equation. We get $ab\in C$ where $b=(\omega\ot\iota)(\co(h)(a'\ot
1))$. This proves that $C$ is an algebra.
\snl
We have that $\co(h)(A\ot 1)$ is a subspace of
$A\ot C$. If we apply $\co\ot\iota$ and multiply with $A\ot 1\ot 1$
we get
$$((\iota\ot\co)\co(h))(A\ot A\ot 1)\subseteq A\ot A\ot C.$$
When we apply a linear functional $\omega$ on the first leg of this
equation, we get
$$\co(a)(A\ot 1)\subseteq A\ot C$$ for all $a\in C$.
\snl
Finally, it is not difficult to see that $h\in C$. Indeed, for any
$a\in A$ with $\varepsilon(a)=1$ we get $h\in C$ from applying
$\varepsilon\ot\iota$ to the expression $\co(h)(a\ot 1)$.
\einspr

We see that the right leg of $\co(C)$ lies again in $C$.
Therefore, $C$ is a {\it left invariant \st subalgebra} of $A$ (see
Definition A.13 in the
 appendix).
\snl
Observe that condition i) in Proposition 3.2, together with the
left invariance of $C$ also implies ii) in that proposition.

\nl
\it The algebra $C$ in the case of an algebraic quantum group \rm
\nl
Now, and further in this section, we will assume that the pair
$(A,\co)$ is an algebraic quantum group. As before, we will use
$\varphi$ for a left integral and $\psi$ for a right integral.
Given the group-like projection $h$, we assume that these two
integrals are normalized so that $\varphi(h)=1$ and $\psi(h)=1$.
\snl
First of all, in this case, we can use the formulas ii) and iii)
in Proposition 3.2 to characterize elements in $C$. This is done
in the following proposition.

\inspr{3.4} Proposition \rm
Let $(A,\co)$ be an algebraic quantum group and let $C$ be
the right leg of a group-like projection $h$ in $A$. Then
$$\align
     C &= \{ a\in A \mid \co(a)(1\ot h) = a\ot h \} \\
       &= \{ a\in A \mid (1\ot a)\co(h) = (S(a)\ot 1)\co(h) \}.
\endalign$$
\snl\bf Proof: \rm
Suppose that $a \in A$ and that $\co(a)(1\ot h)=a\ot h$. We saw in the
proof of Proposition 3.2 that it then follows that also $(1\ot a)\co(h) =
(S(a) \ot 1)\co(h)$. Apply the right invariant functional $\psi$ to
get $\psi(h) a \in C$. As we know that $\psi(h) \neq 0$, we will have
$a\in C$. Then this proposition follows from Proposition 3.2.
\einspr

It is not clear if this result is true in general (that is if we do
not assume the existence of integrals).

\nl
As an application of the above result, let us look again at the case of two
group-like projections $h$ and $h'$ with $h\leq h'$.
\snl
We have seen already in Section 1 that 
$$\co(h')(1\ot h) = h' \ot h \qquad\qquad \text{and}\qquad\qquad
         \co(h')(h\ot 1) = h' \ot h.$$
This is a simple consequence of the definitions and the formula 
$h=h'h$. 
\snl
Now denote by $C_h$ and by $C_{h'}$ the right legs of $\co(h)$ and
$\co(h')$ respectively.

\inspr{3.5} Proposition \rm
With the notation and conditions as above, we have
$C_{h'}\subseteq C_h$.

\snl\bf Proof: \rm
Because $\co(h')(1\ot h)=h'\ot h$, it follows from Proposition 3.4
that $h'\in C_h$. By the left invariance of $C_h$, it follows that
also the right leg of $\co(h')$ sits in $C_h$. This proves the
result.
\einspr

We see that, increasing projections give rise to decreasing
algebras. Compare with the result we found in Proposition 1.9 where we saw that 
a similar sitation happened with the Fourier transform of group-like projections.

\nl
Now again, consider the case of a single group-like projection. We prove some 
more properties about $C$ that are easy consequences
of the properties of $h$ and the definition of $C$.

\inspr{3.6} Proposition \rm
Let $(A,\co)$ be an algebraic quantum group and as before, let $C$ be
the right leg of $\Delta(h)$ for a group-like projection $h$ in $A$. Then $S^2(C)=C$ and
also $\sigma(C)=C$. In fact, $\sigma(a)=S^2(a)$ for all $a\in C$. We
also have $C\delta=\delta C=C$.
\snl\bf Proof: \rm
We have seen (in Proposition 1.6 and Proposition 1.7) that
$\sigma(h)=h$ and $S(h)=h$, so also $S^2(h)=h$. Then
$$\align
    (S^2 \ot S^2)\co(h) &= \co(h) \\
    (S^2 \ot \sigma)\co(h) &= \co(h).
\endalign$$
We use that $\co(\sigma(a))=(S^2 \ot \sigma)\co(a)$ for all $a\in
A$ to get the last formula (see Proposition A.5 in the appendix). 
It follows that indeed $S^2(C)=C$ and
$\sigma(C)=C$. Moreover, combining the two formulas we get
$$ (\iota \ot S^2)\co(h) = (\iota\ot\sigma)\co(h)$$
and this will imply that $\sigma$ and $S^2$ coincide on $C$.
\snl
To prove the second statement, start from the equality
$h\delta=h$, proven in Proposition 1.7 and apply $\co$ to get
$\co(h)(\delta\ot\delta)=\co(h)$. It follows that $C\delta=C$ and
by taking adjoints also that $\delta C=C$.
\einspr

From one of the formulas in Proposition A.9 in the appendix, it follows that for any algebraic quantum group $(A,\Delta)$ with a unimodular dual (i.e.\ when left and right integrals on the dual coincide), we have that $\sigma=S^2$ on $A$. In particular, because a compact quantum group is unimodular, we have this property on a discrete quantum group. In the previous proposition, we see that also $\sigma=S^2$ on this subalgebra $C$. This is not unexpected as we will see later that $C$ shares also other properties with discrete quantum groups.
\snl  
The last formula in the previous proposition says that $\delta$
and $\delta^{-1}$ are also multipliers of $C$.

\nl
In general, we will not have that $C$ is left invariant by the
antipode $S$. We have that $S(C)$ is the left leg of $\co(h)$. It
will also be a \st subalgebra of $A$, now a right invariant one.
We will only have $S(C)=C$ if the two legs of $\co(h)$ are the
same. This is one of the special cases we have considered in
Section 1 (cf.\ Proposition 1.10), but in its dual form. 
We will show that the pair
$(C,\co)$ is a discrete quantum group if $S(C)=C$ (see Theorem 3.8
below). Then we know that we must have that the \st algebra $C$ is
a direct sum of matrix algebras. Now, it turns out that this is
already the case, also when we do not assume that $S(C)=C$. We obtain this
in the next proposition.

\inspr{3.7} Proposition \rm
Let $(A,\co)$ be an algebraic quantum group and $h$ a group-like
projection. Let $C$ be the right leg of $\co(h)$. Then $C$ is a direct
sum of matrix algebras.
\snl\bf Proof: \rm
Take $a\in A$ and write $\co(h)(1 \ot a) = \sum a_i \ot b_i$. Let $V$
be the vector space spanned by the elements $b_i$. We see that
$Ca\subseteq V$. In particular, $Ca$ is finite-dimensional for all
$a\in A$. Similarly, $aC$ is finite-dimensional for all $a$. Then also
$CaC$ will be finite-dimensional for all $a\in A$. So, $CaC$ will be a
finite-dimensional two-sided \st ideal of $C$ for all self-adjoint
elements $a\in C$. Because we are working with operator algebras
(i.e.\ \st algebras with enough positive linear functionals), such a
two-sided \st ideal must be in itself a finite-dimensional operator
algebra, so a direct sum of matrix algebras.
\snl
Now, take any $a\in C$. We know that
$\psi(h)a=(\psi\ot\iota)(\co(h)(1\ot a))$ when $\psi$ is a right integral. 
We also know that there is
an element $e\in A$ such that $(e\ot 1)\co(h)(1\ot a)=\co(h)(1\ot
a)$ (see Proposition A.7 in the appendix). If we apply $\psi\ot\iota$, we see
that $a\in Ca$. Taking adjoints,
 also
$a\in aC$. Combining, we get $a\in CaC$. So any element of $C$ belongs
to such a sum of matrix algebras in $C$.
\snl
Hence $C$ will be a direct sum of matrix algebras.
\einspr

We have seen already that $ah=\varepsilon(a)h$ for all $a\in C$ and so 
$\Bbb C h$ is a one-dimensional component of $C$. The complement of
this component in this direct sum decomposition is precisely the
kernel of $\varepsilon$ on $C$.
\snl
Observe that not very much is needed about $h$, except for the fact that $C$ is a $^*$-algebra. Indeed, in [L-VD2], we will obtain a similar result as above for any left
invariant \st subalgebra with a non-trivial relative commutant in the dual algebra
$B$.

\nl
Let us now combine the result of Proposition A.15 in the appendix (exitence of local
units in $C$) with the one of Proposition 3.7. Let elements $a_1, a_2,
\dots, a_n$ in $A$ be given and let $e\in C$ so that $a_ie=a_i$
for all $i$. Because of the structure of $C$ (being a direct sum
of matrix algebras), we can find a projection $f$ in $C$ such that
$ef=e$. This element will still satifsfy $a_if=a_i$ for all $i$.
It is not difficult to show that in fact, we have two-sided local
units for $A$ that are projections (in $C$). This strengthens the
result in Proposition A.15.
\nl
\it The comultiplication $\co$ on $C$ in the case when $S(C)=C$
\rm
\nl
As mentioned before, the algebra $C$ should be considered as the dual object to $Ah$, or more precisely, as the inverse Fourier transform of $Bk$ where $k$ is the Fourier transform of $h$. This means that the condition $S(C)=C$ should be viewed as dual to the condition $S(Ah)=Ah$, or again more precisely, as $S(Bk)=Bk$. This is the case when $Bk=kB$ since $S(k)=k$. This means that $k$ is central. Therefore it is quite natural to start with this case, just as in Section 1 where we first considered the case where $h$ is central.
\snl
In the special case where  $S(C)=C$, we have that $C$ is left
invariant by $\co$ and that $(C,\co)$ is actually a discrete
quantum group. This is the content of what follows.

\inspr{3.8} Theorem \rm
Let $(A,\co)$ be an algebraic quantum group and assume that $h$ is
a group-like projection such that the left and the right legs of
$\co(h)$ coincide. Then $(C,\co)$ is a discrete quantum group. The
element $h$ is the co-integral.
\snl\bf Proof: \rm
We have $\co(C)(A\ot 1)\subseteq A\ot C$ and by the symmetry also
$\co(C)(1 \ot A)\subseteq C\ot A$. This will give $\co(C)(C\ot
1)\subseteq C\ot C$ and $\co(C)(1\ot C) \subseteq C\ot C$. So, the
maps $T_1$ and $T_2$, defined by $T_1(a\ot b)=\co(a)(1\ot b)$ and
$T_2(a\ot b)=(a\ot 1)\co(b)$ will go from $C\ot C$ to $C\ot C$. Of
course they are injective on $C\ot C$ as they are already injective on
$A\ot A$. The inverse of the map $T_1$ on $A\ot A$ is given by
$$T_1^{-1}(a\ot b)=(\iota\ot S)((1\ot S^{-1}(b))\co(a))$$
and it follows from the fact that $C$ is invariant under $S$, that also
$T_1^{-1}$ will map $C\ot C$ into itself. Similarly for
$T_2$. Therefore $(C,\co)$ is a multiplier Hopf $^*$-algebra.
\snl
Clearly $h$ is a co-integral as $ah=\varepsilon(a) h$ for all
$a\in C$ (and also $ha=\varepsilon(a)h$ by taking adjoints). It
follows that $C$ is a discrete quantum group (Definition 5.2 in [VD3], see also Definition A.8 in the appendix).
\einspr

So, this theorem is a dual version of Theorem 2.2 in the previous
section. In [L-VD2], where we study all these objects in duality, 
we will come back to these two special cases and see how they are dual one to the
other.

\nl
Let us now illustrate this case with two examples.

\inspr{3.9} Examples \rm i) First let $K$ be a finite subgroup of a
(discrete) group $G$. Let $A$ be the group algebra and let
$h=\frac{1}{n}\sum_{p\in K} \lambda_p$ (see Example 1.2.ii). Of
course, as $A$ is co-abelian, $\co(h)$ is symmetric. Here, the
algebra $C$ is nothing else but the group algebra of $K$.
\snl
ii) Let $K$ be an open compact subgroup of a locally compact group
$G$ and let $A$ be the algebra of polynomial functions on $G$ as in Example 1.3.i. Assume that $K$ is normal and let $h=\chi_K$, the characteristic
function of $K$. The right leg of
$\co(h)$ consists of continuous functions with compact support and
constant on right cosets. Because $K$ is assumed to be a normal
subgroup, the right cosets are the same as the left cosets and so we
get that the two legs of $\co(h)$ are the same. The discrete quantum
group we get here is precisely the one associated to the discrete
group $G/K$. The algebra $C$ is now the algebra of complex functions
with finite support on $G/K$.
\einspr

\nl
\it The \st algebra $C_1=C\cap S(C)$ \rm
\nl
Now we will see what is still possible in the general case, that
is when the algebra $C$ is not invariant under the antipode $S$.
It is more or less obvious that we have to consider the intersection
$C \cap S(C)$. This is the step dual to taking $hAh$ in the previous section. Indeed, $hAh=Ah\cap hA$ and $hA=S(Ah)$.

\iinspr{3.10} Proposition \rm
Let $C_1=C \cap S(C)$. Then $C_1$ is a \st subalgebra of
$A$, invariant under the antipode. We still have $AC_1=A$.

\bf\snl Proof: \rm
We know that $C$ and $S(C)$ are \st subalgebras. Therefore, the
intersection is also a \st subalgebra. Because $S^2(C)=C$, we have
$S(C_1)=C_1$.
\snl
We claim that elements $x$ of the form $(\psi\ot\iota)(\co(a)(b\ot
1))$ with $a,b\in C$ belong to $C_1$. Indeed, by the left invariance 
of $C$ we get $x\in C$ because $a\in C$. We also know
that $x=S(y)$ where $y=(\psi\ot\iota)((a\ot 1)\co(b))$ (see
Proposition A.4 in the appendix) and
 so
$y\in C$ because $b\in C$. Hence $x\in C\cap S(C)$ if both $a,b\in
C$.
\snl
Next, take any finite number of elements $a_1, a_2, \dots, a_n$ in
$A$. By Proposition A.15 we have an element $p\in C$ such that
$$(1\ot a_i)\co(h)(p\ot 1)=(1\ot a_i)\co(h)$$
for all $i$. If we apply $\psi$ we get again $a_ie=a_i$ for all $i$
with $e=(\psi\ot\iota)\co(h)(p\ot 1)$. By the result we just proved,
$e\in C_1$. So we find local units for $A$ in $C_1$. Then $AC_1=A$.
\einspr

We see from this proposition that $C_1$ is still imbedded non-degenerately in $A$ (as $AC_1=A$ and by taking adjoints, also $C_1 A=A$). Again, we can consider the multiplier algebra $M(C_1)$ as sitting inside $M(A)$ (see a remark in the appendix, following Proposition A.15). 
\snl
We have $h\in C_1$ because
$h\in C$ and $S(h)=h$. Clearly $S^2(C_1)=C_1$ because
this holds for $C$. And as $\sigma=S^2$ on $C$, the same is true
on the smaller algebra $C_1$. Finally, we still have
$C_1\delta=\delta C_1=C_1$ because this holds for $C$ and because
$S(\delta)=\delta^{-1}$. So, $\delta$ and $\delta^{-1}$ are still
multipliers of $C_1$.
\snl
Of course, also $C_1$ is a direct sum of matrix algebras.
\nl
Now, we would like to consider the comultiplication on $C_1$. We saw in
Theorem 3.8 that $\co$ leaves $C$ invariant if $S(C)=C$. We can
not expect that $\co$ leaves $C_1$ invariant. Just as in the
compact case (Section 2), we have to use a projection map to get
down to $C_1$. As in that case, also here we will no longer have a
comultiplication in the usual sense as it will not be a
homomorphism anymore. It will be a positive linear map, still
satisfying coassociativity (as was the 
case in Proposition 2.3).
\snl
We will use the projection maps $E:A\to C$ and $E':A\to S(C)$ that
we will introduce in the following proposition. Recall that we use
left and right integrals $\varphi$ and $\psi$ resp.\ that are
normalized w.r.t.\ our group-like projection $h$. So
$\varphi(h)=\psi(h)=1$. Remember that $\psi(\,\cdot\,h)=\varphi(\,\cdot\,h)$ because
$h\delta=h$ and that we use $\varphi_0$ for this functional. 

\iinspr{3.11} Proposition \rm Define linear maps
$E,E'$ from $A$ to itself by
$$\align
     E(a) &= (\iota\ot\varphi)(\co(a)(1\ot h))=(\iota\ot\varphi_0)\co(a), \\
     E'(a) &= (\psi\ot\iota)(\co(a)(h\ot 1))=(\varphi_0\ot\iota)\co(a).
\endalign$$
Then $E$ and $E'$ are positive, faithful conditional expectations from
$A$ onto $C$ and $S(C)$ respectively. The antipode converts one to the
other and they commute with each other.

\bf\snl Proof: \rm
First observe that the linear map $E$ is positive because
$$E(a^*a)=(\iota\ot \varphi)((1\ot h)\co(a^*a)(1\ot h))$$
where we have used that $h^2=h=\sigma(h)$. To show that $E$ is
faithful, take $a\in A$ and assume that $E(a^*a)=0$. It follows
from the fact that $\varphi$ is faithful, that then $\co(a)(1\ot
h)=0$. As $h\neq 0$, we must have $a=0$ (by the very definition of a multiplier Hopf algebra). So, $E$ is a faithful
positive map.
\snl
If $a\in C$, then $\co(a)(1\ot h)=a\ot h$ and we see that $E(a)=a$
because $\varphi(h)$ is assumed to be $1$. If $a\in C$ and $b\in A$ we
have
$$\align
    E(ba) &= (\iota\ot\varphi)(\co(b)\co(a)(1\ot h)) \\
          &= (\iota\ot\varphi)(\co(b)(a\ot h)) \\
          &= E(b)a.
\endalign$$
Similarly, or by taking adjoints, we get also $E(ab)=aE(b)$ if $a\in
C$ and $b\in A$.
\snl
Finally, we need to show that $E(A)\subseteq C$. We know from
Proposition A.4 in the appendix that
 $E(a)=S(x)$ where $x=(\iota\ot\varphi)(\co(h)(1\ot a))$ (using that
$\sigma(h)=h$). So $x$ belongs to
 the left leg of $\co(h)$ and consequently $S(x)$ sits in the right leg
of $\co(h)$. We find $E(a)\in C$
 for all $a\in
A$.
\snl
A similar argument will work for $E'$.
\snl
Now, for any $a\in A$ we have
$$\align
    E(S(a))
       &= (\iota\ot\varphi)(\co(S(a))(1\ot h)) \\
       &= (\varphi\ot\iota)(((S\ot S)\co(a))(h\ot 1)) \\
       &= (\varphi\ot\iota)((S\ot S)((h\ot 1)\co(a))) \\
       &= S((\psi\ot\iota)((h\ot 1)\co(a))) \\
       &= S(E'(a)).
\endalign$$
Similarly we can get $E'(S(a))=S(E(a))$. It is also possible to obtain
this from the previous formula, using that $S^2$ commutes with both
$E$ and $E'$.
\snl
Finally, for any $a\in A$, we have
$$\align
    E'(E(a))
       &= E'((\iota\ot\varphi)(\co(a)(1\ot h))) \\
&= (\psi\ot\iota\ot\varphi)((\co\ot\iota)\co(a)(h\ot 1\ot h))
\endalign$$
and by the coassociativity of $\co$ this turns out to be the same as
$E(E'(a))$. This proves the proposition.
\einspr

It follows that the composition $E'E$ will be a positive, faithful conditional
expectation from $A$ onto $C_1$. We can write
$E'E(a)=(\varphi_0\ot\iota\ot\varphi_0)\co^{(2)}(a)$
for all $a\in A$.
\nl
Before we continue, let us make the following observation. One can verify that the maps $E$ and $E'$ are adjoint to the maps from $B$ to $B$ given by $b\mapsto bk$ and $b\mapsto kb$ respectively. This explains part of the results of the above proposition. The composition $EE'$ is the adjoint of the map $b\mapsto kbk$ on $B$. 
\snl
Now, we are almost ready to define a positive
comultiplication $\co_1$ on $C_1$. The following lemma will be provide the last formula to make this possible.

\iinspr{3.12} Lemma \rm
For any $a\in A$, we get
$$(E\ot\iota)\co(a)=(\iota\ot E')\co(a).$$

\snl\bf Proof: \rm
We have by the definitions of $E$ and $E'$
$$\align
    (E\ot\iota)\co(a) &=
    (\iota\ot\varphi_0\ot\iota)((\co\ot\iota)\co(a)),\\
    (\iota\ot\ E')\co(a) &=
    (\iota\ot\varphi_0\ot\iota)((\iota\ot\co)\co(a))
\endalign$$
and the result follows from the coassociativity of $\co$.
\einspr

Before we continue, we must make a remark about the formula in the
lemma. Indeed, we are applying the maps $E\ot\iota$ and $\iota\ot
E'$ to an element in $M(A\ot A)$ and this is not completely
obvious. One solution to this problem is to give a meaning to the
equation in the lemma by multiplying at the right places with
elements of $A$ (or better, with elements in $C_1$). Another, more
elegant solution is to extend these maps to $M(A\ot A)$. We know
how to do this for non-degenerate homomorphisms. In a similar way,
this can be done for these conditional expectations. Take e.g.\ an
element $m\in M(A)$ and define a multiplier $E(m)$ in $M(C)$ by
$$\align E(m)c&=E(mc) \\
         cE(m)&=E(cm)
\endalign$$
when $c\in C$. It is easy to verify that $E(m)$ is a well-defined
element in $M(C)$. Also $E$, as a map from $M(A)$ to $M(C)$,
satisfies the expected properties (like $E(1)=1$). In a similar
way, one can extend the maps $E'$, $E\ot\iota$ and $\iota\ot E'$.
\nl
\it The comultiplication $\co_1$ on $C_1$ \rm
\nl
The following is an easy consequence of Lemma 3.12.

\iinspr{3.13} Proposition \rm
Define $\co_1$ on $A$ by
$$\co_1(a)=(E\ot\iota)\co(a)=(\iota\ot E')\co(a).$$
Then
$\co_1(C_1)(1\ot C_1) \subseteq C_1 \ot C_1$ and 
$\co_1(C_1)(C_1\ot 1) \subseteq C_1 \ot C_1$.
We also have that $\co_1$ is coassociative on $C_1$.

\snl\bf Proof: \rm
We know that $\co(a)(A\ot 1)\subseteq A\ot C$ when $a\in C$. Because
$\co_1=(\iota\ot E')\co$ and $E'(C)\subseteq C_1$ we also get
$\co_1(a)(A\ot 1)\subseteq A\ot C_1$ for all $a\in C$. Similarly,
we know that $\co(a)(1\ot A)\subseteq S(C)\ot A$ when $a\in S(C)$. And now,
because also
$\co_1=(E\ot\iota)\co$ and $E(S(C))\subseteq C_1$ we get
$\co_1(a)(1\ot A)\subseteq C_1\ot A$ for all $a\in S(C)$.
So
$$\co_1(C_1)(1\ot A) \subseteq C_1 \ot A \qquad\quad \text{and}\qquad\quad
    \co_1(C_1)(A\ot 1) \subseteq A\ot C_1.$$
Then also
$$\co_1(C_1)(1\ot C_1) \subseteq C_1 \ot C_1 \qquad\quad \text{and}\qquad\quad
    \co_1(C_1)(C_1\ot 1) \subseteq C_1 \ot C_1.$$
To prove coassociativity we use that
$\co_1(a)=(\iota\ot\varphi_0\ot\iota)\co^{(2)}(a)$
so that
$$(\iota\ot
\co_1)\co_1(a)=(\iota\ot\varphi_0\ot\iota\ot\varphi_0\ot\iota)
\co^{(4)}(a).$$
We get the same right hand side if we calculate $(\co_1 \ot
\iota)\co_1(a)$.
\einspr

This map $\co_1$ will not be a homomorphism, but it will again be
positive. The reason for this is that $\co$ is a \st
homomorphism and that the maps $E,E'$ are positive.
It also follows from the expression
$$\co_1(a)=(\iota\ot\varphi_0\ot\iota)\co^{(2)}(a)$$
and the positivity of $\varphi_0$. Also here, as in the case of
$\co_0$, defined in Proposition 2.3, it
 can be shown that $\co_1$ is completely positive, but again this fact
is of no interest for this paper.

\nl
We have a {\it counit} for this new comultiplication on $C_1$. It is
simply the restriction of the original comultiplication
$\varepsilon$ to $C_1$. We obtain a \st homomorphism from $C_1$ to
$\Bbb C$ still satisfying the axioms for a counit on the new
comultiplication. Indeed, when $a\in C_1$ we have
$$\align (\varepsilon\ot\iota)\co_1(a)
     &= (\varepsilon\ot\iota)(\iota\ot E)\co(a) \\
     &= E((\varepsilon\ot\iota)\co(a)) = E(a) = a
\endalign$$
and similarly $(\iota\ot\varepsilon)\co(a)=a$.
\nl
\it The integrals on the pair $(C_1,\co_1)$ and the antipode. \rm
\nl
First, we look for the left and right integrals on the algebra
$C_1$, considered with the  positive comultiplication
$\co_1$. These present no problem because of the following result.

\iinspr{3.14} Lemma \rm
The left and right integrals $\varphi$ and $\psi$ on $A$ are both invariant
under $E$ and $E'$.
\snl\bf Proof: \rm
We have
$$\align
\psi(E(a)) &= (\psi\ot\varphi_0)\co(a) = \psi(a)\varphi_0(1) = \psi(a)\\
\varphi(E'(a)) &= (\varphi_0\ot\varphi)\co(a) = \varphi_0(1)\varphi(a)
         = \varphi(a).
\endalign$$
On the other hand, because $\varphi_0(\delta)=\varphi(\delta
h)=\varphi(h)=1$ and also $\varphi_0(\delta^{-1})=1$, we have
$$\align
\varphi(E(a)) &= (\varphi\ot\varphi_0)\co(a) =
\varphi(a)\varphi_0(\delta) = \varphi(a)\\
\psi(E'(a)) &= (\varphi_0\ot\psi)\co(a) = \varphi_0(\delta^{-1})\psi(a)
         = \psi(a).
\endalign$$
\einspr

It follows easily from this that the restrictions of $\varphi$ and
$\psi$ to $C_1$ are still left, resp.\ right invariant for the new
comultiplication. Indeed, let e.g. $a\in C_1$. Then
$$ (\varphi\ot\iota)\co_1(a)=(\varphi\circ E\ot\iota)\co(a)
         =(\varphi\ot\iota)\co(a)=\varphi(a)1.$$
Observe also that the automorphisms $\sigma$ and $\sigma'$ leave
$C_1$ globally invariant and so they will still be
automorphisms for the restrictions.
\nl
Now we look at the {\it antipode} and its relation with $\co_1$ and with
the integrals.
\snl
We have seen already that $S$ leaves $C_1$ invariant. Also
$E'(S(a))=S(E(a))$ for all $a\in A$. Then we get, for all $a\in A$,
$$\align \co_1(S(a)) &= (\iota\ot E')\co(S(a)) \\
                     &= (\iota\ot E')(S\ot S)\co^{\text{op}}(a) \\
                     &= (S\ot S)(\iota\ot E)\co^{\text{op}}(a) \\
                     &= (S\ot S)\co_1^{\text{op}}(a)
\endalign$$
as expected. However, we also know that this property does not
characterize the antipode completely. We must, as we did also in
Section 2, look at the relation with the integrals.
\snl
We have the following result.

\iinspr{3.15} Proposition \rm For all $a,b\in C_1$  we have
$$S((\iota\ot\varphi)(\co_1(a)(1\ot b)))=(\iota\ot\varphi)((1\ot
a)\co_1(b)).$$

\snl\bf Proof: \rm
Let $a,b\in C_1$. Then
$$\align \co_1(a)(1\ot b) &= (\iota\ot E')(\co(a))(1\ot b) \\
                          &= (\iota\ot E')(\co(a)(1\ot b))
\endalign$$
where we use that $b\in S(C)$ and therefore $E'(xb)=E'(x)b$ for
all $x\in A$.
\snl
Now, apply the left integral $\varphi$ and use that $\varphi$ is
invariant for $E'$ to get
$$(\iota\ot\varphi)(\co_1(a)(1\ot b))=(\iota\ot\varphi)(\co(a)(1\ot
b)).$$
Similarly we have
$$(\iota\ot\varphi)((1\ot a)\co_1(b))=(\iota\ot\varphi)((1\ot
a)\co(b)).$$ Then the result follows because we already have
$$S((\iota\ot\varphi)(\co(a)(1\ot b)))=(\iota\ot\varphi)((1\ot
a)\co(b))$$
for all $a,b\in A$.
\einspr

Also here it is possible to prove the following result. It is the
dual version of Proposition 2.6 of the previous section.

\iinspr{3.16} Proposition \rm
We have
$$C_1=\text{span}\{(\iota\ot\varphi)(\co_1(a)(1\ot b)) \mid a,b\in
C_1 \}.$$
\einspr

Again, the result follows essentially from the faithfulness of
$\varphi$ and the existence of the
 counit (cf.\ the proof of Proposition 2.6).
\snl
If we combine the previous results, we arrive at the following theorem, which is a dual version of Theorem 2.7 in the previous section.

\iinspr{3.17} Theorem \rm
Let $(A,\Delta)$ be an algebraic quantum group and $h$ a group-like projection in $A$. Let $C_1=C\cap S(C))$ where $C$ is the right leg of $\Delta(h)$ and let $\Delta_1$ be the coproduct defined on $C_1$ as in Proposition 3.13. Then $(C_1,\Delta_1)$ is an algebraic quantum hypergroup of discrete type (as introduced in [De-VD1]).
\einspr

It is of discete type because there is a co-integral, namely $h$ itself. In our second paper on this subject [L-VD2], we will see $(C_1,\Delta_1)$ is the dual version of the compact quantum hypergroup $(A_0,\Delta_0)$ from Theorem 2.7. More precisely, it is in duality with $(B_0,\Delta_0)$, defined in a similar way as $A_0$.  
\nl

Let us now finish this section with the basic examples.

\iinspr{3.18} Examples \rm 
i) Let $K$ be a finite subgroup of a (discrete) group $G$. Consider the group algebra of $G$ and let $h=\frac1n\sum_{p\in K}\lambda_p$ as in Example 1.2.ii (where $n$ is the number of elements in $K$). Because $\Delta$ is coabelian, the left and right legs of $\Delta(h)$ coincide so that $C_1=C$. Then also $\Delta_1=\Delta$ and we are in the situation of Example 3.9.i. The pair $(C,\Delta)$ is nothing else but the group algebra of $K$. A similar situation occurs if we take a compact open subgroup of a locally compact group $G$ and $h$ as in Example 1.3.ii. Now, we will get the $^*$-subalgebra of $C_r^*(K)$ consisting of polynomial functions on $K$. This is a discrete quantum group.
\snl
ii) Again, let $K$ be a finite subgroup of a group $G$, but consider the algebra $A$ of functions with finite support on $G$ and pointwise operations. Take for $h$ the characteristic function of $K$ as in Example 1.2.i. The right leg $C$ of $\Delta(h)$ is the $^*$-algebra of functions $f$ with finite support on $G$ and such that $f(pk)=f(p)$ for all $p\in G$ and $k\in K$. Similarly, $C_1$ will be the $^*$-algebra of functions $f$ with finite support on $G$ and such that $f(kpk')=f(p)$ for all $p\in G$ and $k,k'\in K$. The conditional expectations $E$ and $E'$ are given by
$$\align E(f)(p)&=\frac1n\sum_{k\in K}f(pk) \\ 
              E'(f)(p)&=\frac1n\sum_{k\in K}f(kp) 
\endalign$$ 
where again, $n$ is the number of elements in $K$, and where $f\in A$ and $p\in G$. The coproduct $\Delta_1$ on $C_1$ satisfies 
$$\Delta_1(f)(p,q)=\frac1n\sum_{k\in K}f(pkq).$$
It is easily verified that the pair $(C_1,\Delta_1)$ is indeed an algebraic quantum hypergroup of discrete type.
\snl
iii) In a completely similar way, we can look at the case of a locally compact group $G$ with a compact open subgroup. Now, the algebra $A$ is the algebra  of polynomial functions and also $h$ is the characteristic function of $K$ (as in Example 1.3.i). For $C$ we get the $^*$-algebra of continuous complex functions with compact support on $G$, constant on right $K$-cosets. For $C_1$ we get the $^*$-algebra of continous functions with compact support on $G$, constant on double $K$-cosets. In the first case, this is the algebra of complex functions with finite support on the qoutient space $G/K$, while in the second case, this is the algebra of complex functions with finite support on the double quotient space $K\backslash G/K$. The sums in the previous example become (normalized) integrals over $K$. So, the coproduct $\Delta_1$ on $C_1$ is given by the formula
$$\Delta_1(f)(p,q)=\int_K f(pkq)\,dk$$
where the normalized Haar measure on $K$ is used. 
\einspr 

We can make a remark, similar as Remark 2.10.iv in the previous section. Depending on the point of view, we can think of the algebras $C$ and $C_1$ either as discrete quantum (hyper) qoutient groups or discrete quantum (hyper) subgroups.

\nl
\nl

\bf 4. Conclusions and further research \rm
\nl
In a {\it forthcoming paper} [L-VD2], we study the objects, 
obtained in this paper, in duality. This will only be possible for
algebraic quantum groups. In Section 1 of this paper, we have already shown that the
Fourier transform $k$ of a group-like projection $h$ is again a
group-like projection in the dual $(\hat A,\hat \co)$ of $(A,\co)$.
Therefore, we can also consider the algebras, associated with the
projection $k$ in the dual: $\hat A k$, $k\hat A k$, $D$ and
$D_1$. This will be done in [L-VD2].

\snl
We will look at the relative position of the
algebras $C$ and $D$ in the Heisenberg algebra. Here, the
Heisenberg algebra is the algebra generated by $A$ and its dual
$\hat A$ with the Heisenberg commutation relations between
elements of $A$ and  $\hat A$. It turns out that $C$ and $D$ are
exactly each others relative commutants:
$$\align
    C &= \{a\in A \mid ab=ba \text{ for all } b\in D \} \\
    D &= \{b\in \hat A \mid ab=ba \text{ for all } a\in C \}.
\endalign$$
This result is of the same type as a well-known result in the
theory of crossed products, cf. [L]. We also prove some result in the
other
 direction.
We show that any left invariant \st subalgebra $C$ of $A$ with a
non-trivial relative commutant $D$, must be the right leg of
$\co(h)$ for a group-like projection $h\in A$. It follows that then
$D$ is the relative commutant of $C$. So, we obtain a one-to-one
correspondence between certain left invariant $^*$-subalgebras and
group-like projections.

\snl
In [L-VD2], we also investigate two types of relations between
the objects obtained in Section 2 (the compact quantum (hyper)
subgroups) and those obtained in Section 3 (the discrete quantum
(hyper) subgroups). One aspect is the duality between those two.
The other one is the Fourier transform that maps one object to its
dual object.

\snl
We will also include more examples in [L-VD2] and illustrate our results
using these examples. Related is also the work on algebraic quantum hypergroups [De-VD2].

\nl
\it Further research \rm

\nl
In this paper, we have essentially only treated the case of an
algebraic quantum group. A few results were also shown to be true for all
multiplier Hopf \st algebras and probably more results are correct in this more general situation. It is one possible direction of future research: What
can still be done for general multiplier Hopf algebras? In [L-VD2] and [L-VD3],
when studying special cases and examples, we also go beyond
the case of algebraic quantum groups. It is an indication that
some of the objects we studied and results we obtained, can also
be studied in the more general setting of locally compact quantum
groups. Remember that already in this introduction, we mentioned
that algebraic quantum groups serve as a laboratory for work on
general locally compact quantum groups.

\snl
One of the related problems is that of characterizing the algebraic
quantum groups among the general locally compact quantum groups. There
is now a solution to this problem in the abelian case (i.e.\ for
$C_0(G)$) and in the dual, coabelian case (i.e.\ for the dual of
$C_0(G)$, the reduced group C\st algebra $C_r^*(G)$ of $G$), see [L-VD1]. The
solution in these two special cases may inspire a possible solution in
the general case. A related question that remains open is the following. Is there
always a group-like projection in an algebraic quantum group? Observe 
that the answer is positive in the abelian and the coabelian case.

\snl
Another related problem is to find some structure theorem for
algebraic quantum groups. What we have in mind is the special case of
a locally compact group $G$ with a normal compact open subgroup
$K$. Then $G$ is a cocycle crossed product of $K$ with $G/K$. Similar
constructions exist in the theory of locally compact quantum groups
(see [V-V]). One can hope that this work could be used to understand
the structure of an algebraic quantum group (more likely even when
there exists a group-like projection).

\snl
Totally disconnected locally compact groups have a local basis of
the identity consisting of compact open subgroups. Taking into account the results in this paper, we propose the following definition of a {\it totally disconnected locally compact quantum group}. It should be an algebraic quantum group $(A,\Delta)$ with the property that there exists a decreasing net $(h_\lambda)_\lambda$ of group-like projections in $A$ such that the associated $^*$-algebras $C_\lambda$, the right legs of $\Delta(h_\lambda)$, increase to $A$. Remark that we have shown in Section 3 that indeed, these algebras will give a increasing net. So, we just require that any element of $A$ belongs to such an algeba $C_\lambda$ for some $\lambda$. It is not so hard to show that in this case, also the associated smaller $^*$-algebras $C_\lambda \cap S(C_\lambda)$ will increase to $A$. In any case, it is
 expected that the development of a theory of totally
disconnected locally compact quantum groups will depend heavily on the results in
this paper.

\snl
Finally, not only compact open quantum subgroups of a locally
compact quantum groups can be studied, but also compact quantum
subgroups, closed quantum subgroups, ... and similar questions can
be asked about such subgroups.

\nl\nl

\bf Appendix \rm
\nl
In this appendix, we will first recall the notion of a multiplier
Hopf $^*$-algebra and invariant integrals. Then, we recall some of the main 
properties and in particular, we look at the concept of the Fourier transform. 
All these things are available in the literature (cf.\ the references given in 
the introduction). We also introduce the notion 
of the legs of a comultiplication as this concept is used in this paper. 
Finally we discuss some properties of left invariant algebras as they also 
appear in this work.

\nl
\it Multiplier Hopf $^*$-algebras and algebraic quantum groups\rm

\nl
Let $A$ be an algebra over $\Bbb C$, possibly without identity,
but always with a non-degenerate product. This means that $b=0$ whenever
either $ab=0$ for all $a$ or $ba=0$ for all $a$. This property is automatic when the algebra has an identity.
The {\it multiplier algebra} $M(A)$ of an algebra $A$ can be
characterized as the largest algebra with identity, containing $A$
as a two-sided ideal with the property that, if $x\in M(A)$ then $x=0$ if either
 $xa=0$ for all $a\in A$ or $ax=0$ for all $a\in A$. If $A$ is a \st algebra, then $M(A)$ is
also a \st algebra. The typical example to have in mind is the
algebra of complex functions with finite support on a set (with
pointwise operations). In that case, the multiplier algebra is the
algebra of all complex functions on this set. For a precise
treatment of the multiplier algebra of an algebra, we refer to
[VD3].

\snl
If $A$ is an algebra with a non-degenerate product, the tensor
product $A\ot A$ of $A$ with itself is again an algebra with a
non-degenerate product. We have natural inclusions
$$A\ot A \subseteq M(A) \ot M(A) \subseteq M(A \ot A).$$
In most cases, these two inclusions are strict (when the algebra
has no identity). A {\it comultiplication} on $A$ is a
homomorphism $\co:A \to M(A \ot A)$ satisfying certain properties.
If $A$ has an identity, it is natural to assume that $\co(1)=1 \ot
1$. If the algebra has no identity, this condition is replaced by
non-degeneracy of the comultiplication. This means that
$$\co(A)(A\ot A)=(A\ot A)\co(A)=A\ot A.$$
Then, $\co$ can be uniquely extended to a homomorphism from $M(A)$
to $M(A\ot A)$ and this extension (still denoted by $\Delta$) will be unital. 
We can consider
the homomorphisms $\co\ot\iota$ and $\iota\ot\co$ from $A\ot A$ to
$M(A\ot A\ot A)$ where $\iota$ is used to denote the identity map
from $A$ to itself. These two homomorphisms are again
non-degenerate and can be extended to $M(A\ot A)$. If again we use the
same symbols for these extensions, it makes sense to require
coassociativity for $\co$ in the form
$(\co\ot\iota)\co=(\iota\ot\co)\co$. If $A$ is a \st algebra, also
$A\ot A$ is a \st algebra and in this case, we require $\co$ to be a \st
homomorphism.

\snl
Here, the motivating example is the algebra $A$ of complex
functions with finite support on a group $G$. In this case, $\co$
is defined on $A$ by $(\co(f))(p,q)=f(pq)$ whenever $f\in A$ and
$p,q\in G$. One has to consider elements in $M(A\ot A)$ as complex
functions on $G\times G$. Coassociativity is an easy consequence
of the associativity of the product in $G$.

\nl
We are now ready to recall the definition of a multiplier Hopf
algebra (cf.\ [VD3]).

\inspr{A.1} Definition  \rm
A {\it multiplier Hopf (\st)algebra} is a (\st)algebra $A$
with a non-degenerate product and a comultiplication $\co$
satisfying certain conditions: It is assumed that the two linear
maps, defined from $A\ot A$ to $M(A\ot A)$ by
$$ a \ot b \mapsto \co(a)(1\ot b) \qquad\qquad a \ot b \mapsto (a\ot 1)\co(b),$$
are injective, map into $A\ot A$ and actually have all of $A\ot A$
as their range. A multiplier Hopf algebra is called {\it regular}
if the opposite comultiplication $\co^{\text{op}}$, obtained by
composing $\co$ with the flip, also satisfies these properties. 
\einspr

In the case of a \st algebra, regularity is automatic.
\snl
Again the motivating example is the algebra of complex functions
with finite support on a group $G$ with the comultiplication $\co$
defined as before. The necessary conditions on $\co$ are an easy
consequence of the fact that the maps
$$(p,q) \mapsto (pq,q) \qquad\qquad  (p,q) \mapsto (p,pq)$$
are bijective from $G\times G$ to itself in the case of a group.

\snl
For a multiplier Hopf algebra $(A,\co)$, we have the existence of
a unique {\it counit} and a unique {\it antipode} as in the following
proposition.

\inspr{A.2} Proposition \rm
There exists a unique homomorphism $\varepsilon$ from $A$ to $\Bbb
C$, called the counit, satisfying
$$(\varepsilon\ot\iota)\co(a)=a \qquad\qquad (\iota\ot\varepsilon)\co(a)=a$$
for all $a\in A$. We give a meaning to these equations by
multiplying, left or right, with elements in $A$. There also
exists a unique anti-homomorphism $S$ from $A$ to $M(A)$, called
the antipode, satisfying
$$ m(S\ot\iota)\co(a)=\varepsilon(a)1 \qquad\qquad m(\iota\ot S)\co(a) 
=\varepsilon(a)1$$
for all $a\in A$, where $m$ stands for the multiplication on $A$,
viewed as a linear map from $M(A)\ot A$ or $A\ot M(A)$ to $A$.
Again one multiplies with an element in $A$ to give a meaning to
these formulas. If $A$ is a regular multiplier Hopf algebra, then
$S$ maps into $A$ and it is bijective. If $A$ is a multiplier Hopf
\st algebra, then $\varepsilon$ is a \st homomorphism while $S$
satisfies $S(a^*)=(S^{-1}(a))^*$ for all $a$.
\einspr

Any Hopf (\st)algebra is a multiplier Hopf (\st)algebra.
Conversely, if $(A,\co)$ is a multiplier Hopf (\st)algebra and if
$A$ has an identity, it is a Hopf (\st)algebra. So we see that the
theory of multiplier Hopf algebras extends in a natural way the
theory of Hopf algebras to the case where the underlying algebras
are not required to have an identity.

\snl
For the theory of multiplier Hopf (\st)algebras, we refer to
[VD3] and [VD-Z3].

\nl
Now, let $(A,\co)$ be a regular multiplier Hopf algebra.

\inspr{A.3} Definition \rm
A linear functional $\varphi$ on $A$ is called left invariant if
$(\iota\ot\varphi)\co(a)=\varphi(a)1$ for all $a\in A$. Again
$\iota$ stands for the identity map and the equation is given a
meaning by multiplying with any $b\in A$ from the left. A {\it
left integral} is a non-zero left invariant linear functional on
$A$. Similarly, a right invariant linear functional on
$A$ is defined and it is called a {\it right integral} if it is non-zero.
\einspr

Such integrals on regular multiplier Hopf algebras are unique (up
to a scalar) and if a left integral $\varphi$ exists, also a right integral
exists (namely $\varphi\circ S$). Integrals are automatically faithful. This 
means (for the
left integral $\varphi$) that $\varphi(ab)=0$ for all $a$ will
imply $b=0$ and similarly when $\varphi(ba)=0$ for all $a$.

\snl
If $(A,\co)$ is a multiplier Hopf \st algebra with a positive left
integral, then there is also a positive right integral. This result is non-trivial. It was first shown in [K-VD] but recently, a simpler proof has been obtained in [DC-VD]. 
We will
use the term {\it algebraic quantum group} for such a multiplier
Hopf \st algebra. In this paper, we mainly deal with algebraic
quantum groups, defined in this sense (although probably, some of the results
will still be true without the assumption of positivity of the integrals).

\snl
We refer to [VD7] and [VD-Z3] for details and we will freely
use results from these basic references in this paper. In
particular, we use $\varphi$ to denote a left integral and we use
$\psi$ for a right integral. We use $\sigma$ for the {\it modular
automorphism} of $\varphi$, satisfying
$\varphi(ab)=\varphi(b\sigma(a))$ for all $a,b\in A$. Similarly
$\sigma'$ is used for the modular automorphism of $\psi$. The {\it
modular element} $\delta$ is a multiplier in $M(A)$, defined and
characterized by $(\varphi\ot\iota)\co(a)=\varphi(a)\delta$ for
all $a\in A$. It is invertible and the inverse satisfies
$(\iota\ot\psi)\co(a)=\psi(a)\delta^{-1}$. It is also characterized by the 
formula $\varphi(S(a))=\varphi(a\delta)$ for all $a\in A$. Finally, there is the
{\it scaling constant} $\nu$ defined by
$\varphi(S^2(a))=\nu\varphi(a)$ where $S$ is the antipode. Recently, it is 
shown in [DC-VD] that this scaling constant is trivial ($\nu=1$) for any multiplier 
Hopf $^*$-algebra with positive integrals. There are however examples with non-trivial 
scaling constant if positivity of the integrals is not assumed.

\snl

There are many formulas relating these various objects and they can be
found in the basic references. Nevertheless, let us recall the most important onces, used in this paper. 
\snl
In the first proposition below, we formulate a well-known relation of the 
antipode with the left and the right integrals, important for this paper (see e.g.\ the proof of Proposition 3.11 in [VD7]).

\inspr{A.4} Proposition \rm
For any $a,b\in A$, when 
$$x=(\iota\ot\varphi)(\Delta(a)(1\ot b)) \qquad \text{then} \qquad S(x)=(\iota\ot\varphi)((1\ot a)\Delta(b)).$$ 
Similarly, when 
$$y=(\psi\otimes\iota)((b\ot\iota)\Delta(a)) \qquad \text{then} \qquad S(y)=(\psi\otimes\iota)(\Delta(b)(a\ot 1)).$$
\einspr

Furthermore, we have the following various formulas.

\inspr{A.5} Proposition \rm
We have $\varepsilon(\delta)=1$, $S(\delta)=\delta^{-1}$ and $\sigma(\delta)=\sigma'(\delta)=\nu^{-1}\delta$. Furthermore $\sigma'(a)=\delta\sigma(a)\delta^{-1}$ for all $a$ and $S\sigma'=\sigma^{-1}S$. All the automorphisms $S^2$, $\sigma$ and $\sigma'$ mutually commute and we have
$$\align \Delta(\sigma(a)) &=(S^2 \ot \sigma) \Delta(a) \\
\Delta(\sigma'(a)) &=(\sigma' \ot S^{-2}) \Delta(a)
\endalign$$
for all $a$ in $A$. In the case of a multiplier Hopf $^*$-algebra, we have $\sigma(a)^*=\sigma^{-1}(a^*)$ for all $a$.
\einspr

All these formulas can be found already in the original paper [VD7]. There is also the following result. It was first proven in [K-VD], Lemma 3.10, but we will give here another (new and simpler) argument.

\inspr{A.6} Proposition \rm For all $a\in A$ we have
$$\Delta(S^2(a)) =(\sigma \ot {\sigma'}^{-1}) \Delta(a).$$
 
\snl\bf Proof: \rm
Apply $\varepsilon$ on the second leg in the formula for $\Delta\circ\sigma$ to get 
$$S^{-2}\sigma(a)=(\iota\ot \varepsilon\circ\sigma)\Delta(a).$$ 
Next apply $\varepsilon$ on the first leg in the formula for $\Delta\circ\sigma'$ to obtain
$$S^2\sigma'(a)=(\varepsilon\circ\sigma' \ot \iota)\Delta(a).$$ 
Because $\varepsilon(\delta)=1$, it follows from $\sigma'(a)=\delta\sigma(a)\delta^{-1}$ that $\varepsilon\circ\sigma=\varepsilon\circ\sigma'$. Therefore 
$$(\iota\ot\varepsilon\circ\sigma \ot \iota)\Delta^{(2)}(a)=(\iota\ot\varepsilon\circ\sigma' \ot \iota)\Delta^{(2)}(a)$$
and if we use the two previous formulas, we get from this that 
$$(S^{-2}\sigma \ot\iota)\Delta(a)=(\iota\ot S^2\sigma')\Delta(a).$$
Because $\Delta(S^2(a))=(S^2 \ot S^2)\Delta(a)$, we get the desired formula.
\einspr

For an algebraic quantum group, we have the existence of {\it local units} as in the following proposition.

\inspr{A.7} Proposition \rm
Let $(A,\Delta)$ be an algebraic quantum group. For any finite number of elements $\{a_1, a_2, \ldots, a_n\}$ in $A$ there exists an element $e\in A$ such that $a_je = a_j$ and $ea_j=a_j$ for all $j$.
\einspr

The proof of this result is found in [Dr-VD-Z] and the result is used a few times in this paper (in the proof of Proposition 3.7 and further in this appendix). 
\snl
In this paper, also the following special cases of algebraic quantum groups play a role.

\inspr{A.8} Definition \rm Let $(A,\Delta)$ be an algebraic quantum group. It is called of {\it compact type} if the algebra $A$ has an identity. It is called of {\it discrete type} if there is a left co-integral.
\einspr

Recall that a left co-integral is a non-zero element $h$ in $A$ such that $ah=\varepsilon(a)h$ for all $a$ in $A$. If a left co-integral exists, there is also a right co-integral $k$ satisfying $ka=\varepsilon(a)k$, namely $k=S(h)$. If co-integrals exist, they are also unique up to a scalar.

\nl
\it The dual of an algebraic quantum group and the Fourier transform \rm

\nl
For an algebraic quantum group $(A,\co)$, one can define a dual object $(\widehat A, \widehat \co)$. As a vector space, $\widehat A=\{\varphi(\,\cdot\,a) \mid a\in A\}$. Because of the existence of $\sigma$, it is also true that $\widehat A=\{\varphi(a\,\cdot\,) \mid a\in A\}$. Moreover, in these two cases, the left integral $\varphi$ can be replaced by the right integral $\psi$ (by applying the antipode). The product on $\widehat A$ is dual to the coproduct on $A$ while the coproduct $\widehat \co$ is dual to the product in $A$. It is again an algebraic quantum group. The dual right integral $\widehat \psi$ is given by the formula $\widehat\psi(\omega)=\varepsilon(a)$ when $\omega=\varphi(\,\cdot\,a)$. Observe that the dual of a compact type is of discrete type and vice versa.

\snl
In this paper, we use $B$ for $\widehat A$ and we equip $B$ with the opposite coproduct $\widehat\co^{\text{op}}$. We will now denote this coproduct on $B$ by $\co$. We get a dual pair of multiplier Hopf algebras $(A,B)$ which is somewhat different than in the paper [Dr-VD]. The pairing is denoted by $\langle a, b \rangle$ for $a\in A$ and $b\in B$. It is one of the consequences of this choice that $\langle S(a), b \rangle=\langle a, S^{-1}(b)\rangle$ for all $a\in A$ and $b\in B$. The dual right integral $\widehat\psi$ becomes the left integral on $B$.  The reason for making this choice is to make notations in accordance with the ones that are common in the general theory of locally compact quantum groups as developed in [K-V1] and [K-V2]. 

\snl
It is important to notice that the Sweedler notation can safely be used, also in the theory of multiplier Hopf algebras. One can write $a_{(1)}\ot a_{(2)}$ for $\Delta(a)$ provided we have a {\it covering} of  the first or the second factor, either through a product, from the left or the right, by an element in the algebra, or through a pairing with an element in the dual algebra $B$. We have used the Sweedler notation e.g.\ in the proofs of Proposition 1.6 and 1.8 of Section 1 (and also further in the paper). We refer to [Dr-VD] and [Dr-VD-Z] for more information about the use of the Sweedler notation and this technique of covering (see also [VD-Z3]).

\snl
All the objects that we have introduced for $A$, like $\delta$, $\sigma$ and $\sigma'$, also exist for $B$. We will use the same symbols in most cases as there will be no confusion. Only for the modular elements, we have to be a little more attentive. Remark that in general, the scaling constant for $B$ is $\nu^{-1}$ when $\nu$ is the scaling constant for $A$. This is not so important for this paper because the algebraic quantum groups that are studied here have a trivial scaling constant (see Proposition 1.7 and following comments).

\snl
We will not only have all the formulas for these objects in $B$ as we have for $A$, we also get some new ones, relating the various objects for $A$ with those of $B$. Many of these results can be found e.g.\ in [Ku]. We collect some of them in the following proposition. The results are related with the proof we have given before of Proposition A.6.

\inspr{A.9} Proposition \rm
For the modular element $\delta$ in $B$ we have $\langle a,\delta\rangle=\varepsilon(\sigma(a))=\varepsilon(\sigma'(a))$ and $\langle a,\delta^{-1}\rangle=\varepsilon(\sigma^{-1}(a))=\varepsilon({\sigma'}^{-1}(a))$
for all $a\in A$. Furthermore, for the modular automorphism on $B$, we get $\langle a, \sigma(b)\rangle=\langle \delta^{-1}S^{-2}(a),b \rangle$ and $\langle a, \sigma^{-1}(b)\rangle=\langle \delta S^2(a),b \rangle$ for $a\in A$ and $b\in B$.

\snl\bf Proof: \rm
Let $a, a'\in A$ and $b=\varphi(\,\cdot\,a)$ so that $\varphi(b)=\varepsilon(a)$. If we let $b'=\langle a',b_{(2)}\rangle b_{(1)}$, then one can calculate that $b'= \varphi(\,\cdot\,a\sigma(a'))$ and so $\varphi(b')=\varepsilon(a)\varepsilon(\sigma(a'))$. Then, from the formula $(\varphi\ot\iota)\Delta(b)=\varphi(b)\delta$, it will follow that $\langle a',\delta\rangle=\varepsilon(\sigma(a'))$. It follows from the relation $\sigma'(a)=\delta\sigma(a)\delta^{-1}$ (see Proposition A.5) and $\varepsilon(\delta)=1$ that $\varepsilon(\sigma(a))=\varepsilon(\sigma'(a))$ for all $a$ and this proves the first part of the proposition.
\snl
We will prove the second part of the proposition first in its dual form. 
We start with the formula $S^{-2}\sigma(a)=(\iota\ot\varepsilon\circ\sigma)\Delta(a)$ (cf.\ the proof of A.6). Pairing with an element $b\in B$ and combining with the first result in the proposition, we get precisely $\langle S^{-2}\sigma(a), b\rangle=\langle a, b\delta \rangle$ for all $b$. Moving $S$ to the other side, we get $\langle\sigma(a),b\rangle = \langle a, S^{-2}(b)\delta\rangle$ for all $b$. If we dualize this formula, we get $\langle a, \sigma(b)\rangle=\langle \delta^{-1}S^{-2}(a),b \rangle$, proving the first part of the desired result. The second part is obtained in a similar way, but it also follows easily from the first formula.
\einspr

We must make a few comments before we continue. First, the formula $\langle a,\delta\rangle=\varepsilon(\sigma(a))$ involves a pairing between an element in $A$ and an element in $M(B)$. This can be given a meaning, but the formula can also be interpreted by considering the second form $\langle a', b\rangle=\langle a, b\delta \rangle$ with $a'=(\iota\ot\varepsilon\circ\sigma)\Delta(a)$. 
\snl
Similarly, we find formulas involving $\sigma'$ (using e.g.\ that $S\sigma'=\sigma^{-1}S$).

\nl
Let us now formulate some results about the so-called Fourier transform in this context. These results are not so essential for this paper, but they are of great importance in the second paper [L-VD2]. 
\snl
The {\it Fourier transform} $F$ is considered here as a map from $A$ to $B$ and it is defined by $F(a)=\varphi(\,\cdot\,a)$. Let us fix the left integral $\varphi$ on $B$ with the formula given before and so $\varphi(b)=\varepsilon(a)$ when $b=F(a)$. If $b=F(a)$, then it is not hard to show that $a=\varphi(S^{-1}(\,\cdot\,)b)$ giving the formula for the inverse transform. Indeed, let $a\in A$ and $b=\varphi(\, \cdot \, a)$. If $x\in A$ and $y\in B$ we find
$$\align\langle x, S^{-1}(y)b\rangle &= \langle x_{(1)},S^{-1}(y)\rangle \varphi(x_{(2)}a) \\
                                     &= \langle a_{(1)}, y \rangle \varphi(x a_{(2)})
\endalign$$
using Proposition A.4. Then, using the definition of $\varphi$ on $B$, we find
$$\varphi(S^{-1}(\,\cdot\,)b) = \varepsilon (a_{(2)}) a_{(1)} = a.$$
In the case of a $^*$-algebra, when we assume a positive left integral $\varphi$ on $A$, we get $\varphi(b^*b)=\varphi(a^*a)$ when $b=\varphi(\,\cdot\,a)$. Indeed, using the result above, we get
$$\align \varphi(b^*b) &= \varphi(S^{-1}(S(b^*))b) = \langle a, S(b^*) \rangle\\
      &= \langle S^{-1}(a),b^*\rangle = \langle a^*,b \rangle^- = \varphi(a^*a).
\endalign$$
This is 'Plancherel's formula'. It shows that the dual left integral on $B$ is positive when the left integral on $A$ is positive.
\snl
Finally, observe that, if $A$ has an identity (i.e.\ when it is of compact type), then we have that the left integral $\varphi$ on $A$ is in $B$. It is clearly a left co-integral in $B$ (and so $B$ is of discrete type). On the other hand, if $A$ is of discrete type and if $h$ is a left co-integral in $A$, we see easily that $\varphi(\,\cdot\, h)$ will be the identity in $B$ (provided $\varphi(h)=1$). Remark that $\varphi(h)$ can not be $0$ as this would imply that $\varphi(ah)=0$ for all $a$ and this contradicts the fact that $\varphi$ is faithful.  Hence, $B$ is of compact type.

\nl
\it The left and right legs of $\co(a)$ for $a\in A$ \rm
\nl
In our paper, we use the notion of the legs of $\co(a)$ for $a$ in
a regular multiplier Hopf algebra. We will discuss this concept
here.
\snl
Again, let $(A,\co)$ be a regular multiplier Hopf algebra.
\snl
We will first say what we precisely mean by 'the right leg'
(respectively 'the left leg') of $\co(a)$ for a single element
$a\in A$ or for $a$ belonging to a subspace $V\subseteq A$. We begin with a lemma that prepares for this definition.

\iinspr{A.10} Lemma \rm
Let $a\in A$ and assume that $C$ is a subspace of $A$. Then, the
following are equivalent:
\smallskip
\item{} i) $\co(a)(a'\ot 1)\in A\ot C$ for all $a'\in A$,
\item{} ii) $(a'\ot 1)\co(a)\in A\ot C$ for all $a'\in A$.
\snl\bf Proof: \rm
Suppose i) holds. Take $a'\in A$ and write $(a'\ot 1)\co(a)=\sum
a_i\ot b_i$. Choose $e\in A$ such that $a_ie=a_i$ for all $i$
(cf.\ Proposition A.7). Then
$$\align (a'\ot 1)\co(a) &= (a' \ot 1)\co(a)(e\ot 1) \\
                   &\subseteq (a'\ot 1)(A\ot C) \subseteq A\ot C.
\endalign$$
Therefore ii) holds. Similarly, ii) implies i) and this proves the
lemma.
\snl

\iinspr{A.11} Definition \rm
The smallest subspace $C$ with the properties in the lemma will be
called {\it the right leg} of $\co(a)$. In a similar way, we can
define the left leg of $\co(a)$ for a single element $a\in A$, as
well as the left and right legs of $\co(V)$ for a subspace $V$ of
$A$.
\einspr

It is immediatley clear that the legs of a self-adjoint element
and of a self-adjoint subspace (in the case of a multiplier Hopf
\st algebra) are again self-adjoint subspaces.
\snl
In the case of an algebraic quantum group, we can prove a bit
more.

\iinspr{A.12} Lemma \rm
Let $(A,\co)$ be an algebraic quantum group and let $B$ be the dual. Let
$a\in A$ and let $C$ be a subspace of $A$. Then condition i) and
ii) in Lemma A.10 are also equivalent with the following condition:
\smallskip
\item{} iii) $\langle a_{(1)}, b \rangle a_{(2)} \in C$ for all
$b\in B$.
\snl\bf Proof: \rm
Assume that conditon i) of Lemma A.10 holds. Take any $a'\in A$ and
apply the right integral $\psi$ on the first leg of $\co(a)(a'\ot
1)$. We see that $(\psi(\,\cdot\,a')\ot\iota)\co(a)\in C$ and so
$\langle a_{(1)} , b \rangle a_{(2)} \in C$ for all $b\in B$. This
proves iii). Conversely, assume iii) and take $a'\in A$. Write
$\co(a)(a'\ot 1)=\sum p_i\ot q_i$ with the $(p_i)$ linearly
independent. By assumption, $\sum \psi(p_i x)q_i\in C$ for all
$x\in A$. Let $\omega$ be a linear functional on $A$ such that
$\omega|{_C}=0$. Then $\sum\psi(p_i x)\omega(q_i)=0$ for all $x\in
A$ and by the faithfulness of $\psi$, we have
$\sum\omega(q_i)p_i=0$. As the $(p_i)$ are assumed to be linearly
independent, we must have $\omega(q_i)=0$ for all $i$. Hence
$q_i\in C$ for all $i$. So $\co(a)(a'\ot 1)\in A\ot C$. Therefore
also iii) implies condition i) of Lemma A.10.
\snl
We see that  condition i) of Lemma A.10 is equivalent with iii)
and so the lemma is proved.
\einspr

So in the case of an algebraic quantum group, the right leg of $\Delta(a)$ can also be characterized as the smallest subspace $C$ of $A$ such that $\langle a_{(1)}, b \rangle a_{(2)} \in C$ for all
$b\in B$.

\nl
\it Left invariant subalgebras \rm

\nl
As before let $(A,\co)$ be a regular multiplier Hopf algebra. Using
the notion of the right leg, we can define the following:

\iinspr{A.13} Definition \rm
Let $C$ be a subalgebra of $A$. Then it is called {\it left
invariant} if the right leg of $\co(C)$ belongs to $C$.
\einspr

We call it {\it left} invariant because
$(\omega\ot\iota)\co(C)\subseteq C$ for appropriate linear functionals
$\omega$ on $A$.
\snl
We now prove some useful properties of such left invariant
subalgebras.

\iinspr{A.14} Proposition \rm
Let $C$ be a (non-trivial) left invariant subalgebra of a regular
multiplier Hopf algebra. Then the product in $C$ is (still)
non-degenerate.
\snl\bf Proof: \rm
Suppose that $b\in C$ and that $ba=0$ for all $a\in C$. Then
$(1\ot b)\co(a)(a'\ot 1)=0$ for all $a\in C$ and $a'\in A$. So,
$(1\ot b)\co(a)=0$ for all $a\in C$. If $C$ is non-trivial, this
already implies that $b=0$ (cf.\ Definition A.1). Similarly on the other side.
\einspr

There is even a stronger result. Also for elements $a\in A$ we
have that $a=0$ if $ac=0$ for all $c\in C$ and similarly on the
other side (the proof is as above). This means that $A$ is a {\it
non-degenerate} $C$-bimodule in the sense of [Dr-VD-Z].

\snl
Because the product in $C$ is non-degenerate, it makes sense to
consider the multiplier algebras $M(C)$ and $M(A\ot C)$. It also
follows that $\co(C)\subseteq M(A\ot C)$ for a left invariant
subalgebra. It is not clear however if, for a general subalgebra
$C$ (with a non-degenerate product) of $A$, this condition will be
sufficient to guarantee that we also have the stronger property,
namely that the right leg of $\co(C)$ is contained in $C$. Indeed,
assume that $C$ is a subalgebra of $A$ with a non-degenerate
product and such that $\co(C)\subseteq M(A\ot C)$. Then it is not
hard to see that the second leg of $\co(C)$ belongs to the set
$C_1$, defined by
$$C_1=\{x\in A \mid xa\in C \text{ and } ax\in C \text{ for all } a\in C
\}.$$ In general, it is not clear if $C_1$ is actually contained
in $C$. We say a little more about this after the next
proposition.

\snl
In the following proposition, we show that we always have local
units in a left invariant subalgebra.

\iinspr{A.15} Proposition \rm
Let $(A,\co)$ be an algebraic quantum group and let $C$ be a (non-trivial) left
invariant \st subalgebra of $A$. Given elements $a_1, a_2, \dots,
a_n$ in $A$, there exists an element $e\in C$ so that $a_ie=a_i$
for all $i$.
\snl\bf Proof: \rm
Take any element $b$ in $C$ such that $\psi(b)\neq 0$ where, as
before, $\psi$ is a positive right integral. This is possible
because $C$ is a (non-trivial) \st subalgebra and $\psi$ is a
faithful positive linear functional. Next choose $p\in A$ so that
$$(1\ot a_i)\co(b)(p\ot 1)=(1\ot a_i)\co(b)$$
for all $i$. Then apply the right invariant integral $\psi$ on the
first leg of this equation. This will give $a_ie=a_i$ for all $i$
with $e=(\psi\ot\iota)(\co(b)(p\ot 1))$. Then $e\in C$ and the
result is proven.
\einspr

A consequence of this result is that $AC=A$. Similarly $CA=A$. In
other words, $A$ is a {\it unital} $C$-bimodule, as defined in
[Dr-VD-Z]. It follows that the imbedding of $C$ in $A$ extends to a \st
homomorphism of $M(C)$ to $M(A)$. This extension will be injective
and so we may view $M(C)$ as sitting inside $M(A)$.

\snl
We now come back to a problem discussed earlier (preceding Proposition A.15). If $C$ is a \st
subalgebra satisfying $\co(C)\subseteq M(A\ot C)$ and if $A$ has
local units in $C$, then $C$ will be left invariant in the sense
of Definition A.13. Indeed, we have seen that the right leg of
$\co(C)$ will consist of elements $x\in A$ satisfying $xC\subseteq
C$. Now we can choose $e\in C$ such that $xe=x$. So $x\in C$.
\snl
In Section 3 of this paper, we work with such left invariant $^*$-subalgebras. Because nowhere else in the literature, this concept is defined properly for multiplier Hopf algebras, we have included such a precise definition here in this appendix. 
\nl\nl

\bf References \rm
\bigskip
\parindent 1 cm

\item{[A]} E.\ Abe: \it Hopf algebras. \rm Cambridge University Press (1977).
\smallskip

\item{[C-V]} Y.A.\ Chapovsky \& L.I.\ Vainerman: {\it Compact
quantum hypergroups}. J.\ Operator Theory 41 (1999), 261--289.
\smallskip


\item{[DC-VD]} K.\ De Commer \& A.\ Van Daele. Multiplier Hopf algebras imbedded in C$^*$-algebraic quantum groups. Preprint K.U.\ Leuven (2006). Arxiv math.OA/0611872.
\smallskip


\item{[De-VD1]} L.\ Delvaux \& A.\ Van Daele: {\it Algebraic quantum hypergroups}. Preprint University of Hasselt and K.U.\  Leuven (2006). Arxiv math.RA/0606466. 
\smallskip

\item{[De-VD2]} L.\ Delvaux \& A.\ Van Daele: {\it Algebraic quantum hypergroups. Examples and special cases}. Preprint University of Hasselt and K.U.\  Leuven (in preparation). 
\smallskip


\item{[Dr-VD]} B.\ Drabant \& A.\ Van Daele:
{\it Pairing and the quantum double of multiplier Hopf algebras}.
Algebras and representation theory 4 (2001), 109-132.
\smallskip

\item{[Dr-VD-Z]} B.\ Drabant, A.\ Van Daele \& Y.\ Zhang:
{\it Actions of Multiplier Hopf Algebras.}
Commun.\ in Alg.\ 27(9) (1999), 4117-4127.
\smallskip


\item{[G-H-J]} F.M.\ Goodman, P.\ de la Harpe \& V.F.R\ Jones: {\it Coxeter graphs and towers of algebras}. MSRI publications 14 (1989) Springer, New York.
\smallskip


\item{[Kr]} A.\ Krieg: {\it Hecke algebras}. Mem.\ Amer.\ Math.\ Soc.\ {\bf 87} (1990) no.435.
\smallskip

\item{[Ku]} J.\ Kustermans: {\it The analytic structure of algebraic
quantum groups.} J.\ of Alg.\  {\bf 259} (2003), 415--450.
\smallskip

\item{[K-V1]} J.\ Kustermans \& S.\ Vaes: \it A simple
definition  for locally compact quantum groups. \rm  C.R. Acad. Sci., Paris,
S\'er. I 328 (10) (1999), 871-876.
\smallskip

\item{[K-V2]} J.\ Kustermans \& S.\ Vaes: \it Locally compact
quantum groups. \rm Ann.\ Sci.\
 Ec.\ Norm.\ Sup.\  33
(2000), 837-934.
\smallskip

\item{[K-V3]} J.\ Kustermans \& S.\ Vaes: \it Locally compact
quantum groups in the von Neumann algebra setting. \rm Math.
Scand. {\bf 92} (2003),68--92.
\smallskip

\item{[K-VD]} J. Kustermans \& A.\ Van Daele:
{\it C$^*$-algebraic quantum groups arising from algebraic quantum
groups.} Int.\ J.\ Math.\ 8 (1997), 1067-1139.
\smallskip


\item{[L]} M.\ B.\  Landstad: {\it Duality for covariant systems.} Trans.
Amer.\ Math.\ Soc.\ 248 (1979) 223-267.
\smallskip

\item{[L-VD1]} M.\ B.\ Landstad \& A.\ Van Daele: {\it Groups with compact open subgroups and multiplier
Hopf \st algebras}. Preprint University of
Trondheim and K.U.\ Leuven (2007). Arxiv math.OA/0701525
\smallskip

\item{[L-VD2]} M.\ B.\ Landstad \& A.\ Van Daele: {\it Compact and discrete subgroups of algebraic quantum groups II}. Preprint University of Trondheim and K.U.\ Leuven (2007).
\smallskip




\item{[S]} M.E.\ Sweedler: \it Hopf algebras. \rm Mathematical Lecture Note
Series. Benjamin (1969).
\smallskip


\item{[V-V]} S.\ Vaes \& L.\ Vainerman: {\it Extensions of locally compact quantum groups and the 
bicrossed product construction}. Adv. in Math. 17 (2003), 1--101 .
\smallskip



\item{[VD1]} A.\ Van Daele: {\it Multiplier Hopf algebras.} Trans.\
Amer.\ Math.\ Soc.\ 342 (1994), 917-932.
\smallskip


\item{[VD2]} A.\ Van Daele: {\it Discrete quantum groups.} J.\ of Alg.\
180 (1996), 431-444.
\smallskip


\item{[VD3]} A.\ Van Daele: {\it An algebraic framework for group duality.}
Adv. in Math.  140 (1998), 323--366.
\smallskip


\item{[VD4]} A.\ Van Daele:
{\it The Haar measure on some locally compact quantum groups.}
Preprint K.U.\ Leuven (2001). Arxiv math.OA/0109004.
\smallskip

\item{[VD5]} A.\ Van Daele: {\it Multiplier Hopf $^*$-algebras with
positive integrals: A laboratory for locally compact quantum groups.} Proceedings of the Meeting of Theoretical Physicists and Mathematicians, Strasbourg, February 21 - 23, 2002., Ed. L. Vainerman, IRMA Lectures on Mathematics and Mathematical Physics, Walter de Gruyter, Berlin, New York (2003), pp. 229--247. 
\smallskip

\item{[VD6]} A.\ Van Daele: {\it Locally compact quantum groups: The von Neumann algebras versus the C$^*$-algebra approach.} Bulletin of  Kerala Mathematics Association. Special issue (2006), 153--177 . 
\smallskip

\item{[VD7]} A.\ Van Daele: {\it Locally compact quantum groups. A von Neumann algebra approach.} Preprint K.U.\ Leuven (2006). Arxiv math.OA/0602212
\smallskip




\item{[VD-Z]} A.\ Van Daele \& Y.\ Zhang:
{\it A survey on multiplier Hopf algebras.}
Proceedings of the conference in Brussels on Hopf algebras. Hopf
Algebras and Quantum Groups, eds. Caenepeel/Van Oystaeyen (2000),
269-309. Marcel Dekker (New York).
\smallskip

\item{[W]} S.L.\ Woronowicz: {\it  Compact Quantum Groups.}
Quantum symmetries/Symm\'{e}tries quantiques.  Proceedings of the
Les Houches summer school 1995, North-Holland, Amsterdam (1998),
845--884.
\smallskip

\end